\documentclass[draft,preprint,12pt,numbers,sort&compress]{elsarticle}
\usepackage{mathrsfs}
\usepackage{amssymb}
\usepackage{amsthm}
\usepackage{mathrsfs}
\usepackage[centertags]{amsmath}
\usepackage{amsfonts}

\usepackage{color}

\input amssym.def
\input amssym.tex

\date{}

\newtheorem{Theorem}{Theorem}[section]

\newtheorem{Lemma}{Lemma}[section]

\newcommand\R{\mbox{\bf R}}

\newcommand\SR{\mbox{\scriptsize\bf R}}

\newcommand{\definition}{{\lower .5ex
  \hbox{$\>\>\stackrel{\triangle}{=}\>\>$} }}

\begin{document}

\textwidth=170  true mm
\textheight=270 true mm
\topmargin=-3.0 true cm
\oddsidemargin=0 true cm
\date{}
\mbox{}
\bigskip

\begin{center}{\Large\bf The Cauchy Problem  for Stochastic Generalized  Benjamin-Ono Equation }\\[1ex]
{ Wei YAN$^\dag$}\\[2ex]
{$^\dag$School of Mathematics and Information Science, Henan Normal University,}\\
{Xinxiang, Henan 453007, P. R. China}\\
{ Email: yanwei19821115@sina.cn}\\[2ex]
{Jianhua Huang}\\[2ex]
{$^*$ College of Science, National University of Defense and Technology,}\\
{ Changsha, P. R. China\quad  410073}\\[2ex]
{Email: jhhuang32@nudt.edu.cn}\\
and

{Boling Guo}\\[2ex]
{$^\ddag$Institute of Applied Physics and Computational Mathematics, Beijing, 100088}\\
{Email: gbl@iapcm.ac.cn}
\end{center}

\bigskip
\bigskip

\noindent{\bf Abstract.}  The current  paper  is devoted to  the Cauchy problem for the stochastic  generalized
Benjamin-Ono equation. By establishing the  bilinear estimate, trilinear estimates in some Bourgain spaces,  we prove that
the Cauchy problem for the stochastic generalized  Benjamin-Ono equation  is
locally well-posed for the initial data $u_{0}(x,\omega)\in L^{2}
(\Omega; H^{s}(\R))$ which is $\mathscr{F}_{0}$ measurable with $s\geq\frac{1}{2}-\frac{\alpha}{4}$ and $\Phi \in L_{2}^{0,s}.$
In particular, when $\alpha=1,$ we prove that it is globally well-posed for the initial data
 $u_{0}(x,\omega)\in L^{2}
(\Omega; H^{1}(\R))$ which is $\mathscr{F}_{0}$ measurable and $\Phi \in L_{2}^{0,1}.$
The key ingredients that we use in this paper are trilinear estimates, It$\hat{o}$ formula and the BDG inequality as well as the stopping time technique.
\bigskip


\bigskip
\noindent {\bf Short Title:} Stochastic Generalized  Benjamin-Ono equation
 equation

\bigskip

\leftskip 0 true cm \rightskip 0 true cm

\newpage{}

\newpage{}

{\large\bf 1. Introduction}
\bigskip

\setcounter{Theorem}{0} \setcounter{Lemma}{0}

\setcounter{section}{1}

In this paper, we consider the following stochastic fractional
Benjamin-Ono type equation
\begin{eqnarray}
\label{1.01}
\begin{cases}
du(t)=[|\partial_{x}|^{\alpha+1}\partial_{x}u(t)-\frac{1}{k}\partial_{x}(u^{k})]dt+\Phi dW(t),\\
u(0)=u_0,
\end{cases}
\end{eqnarray}
where
$W(t)=\frac{\partial B}{\partial x}=\sum_{j=1}^{\infty}\beta_{j}e_{j},$ $e_{j}$
is an orthonormal basis of $L^{2}(\R)$ and $(\beta_{j})_{j\in N}$
is a sequence of mutually independent real Brownian motions in a fixed probability
 space and is a  Wiener process on $L^{2}(\R).$
In fact, (\ref{1.01})  is equivalent to the following equations:
\begin{eqnarray}
\label{1.02}
\begin{cases}
\frac{du(t)}{dt}=[|\partial_{x}|^{\alpha+1}\partial_{x}u(t)-
\frac{1}{k}\partial_{x}(u^{k})]+\Phi \frac{dW(t)}{dt},\\
u(0)=u_0.
\end{cases}
\end{eqnarray}
(\ref{1.02})  is considered as the Benjamin-Ono type equation
\begin{eqnarray}
\label{1.03}
\begin{cases}
\frac{dv(t)}{dt}=[|\partial_{x}|^{\alpha+1}\partial_{x}u(t)
-\frac{1}{k}\partial_{x}(u^{k})],\\
u(0)=u_0.
\end{cases}
\end{eqnarray}
forced by a random term $\Phi \frac{dw(t)}{dt}$.

When $\alpha=1$ and $k=2,$  (\ref{1.03})  reduces to the KdV equation  which has been
investigated by many authors, we refer the readers to
\cite{Bourgain,KPV1991, KPVIUMJ, CKST, CKSTEJDE,
CKSTMRL, CKSTJAMS, GD,GJMPA,KPV1993,KPV1996, KPV2001,Kis}.
The result of \cite{KPV1996} and \cite{KPV2001} implies that $s=-\frac{3}{4}$
 is the critical well-posedness indices of the Cauchy problem for the KdV equation.
 Guo \cite{GJMPA} and Kishimoto \cite{Kis} almost proved that the KdV equation is globally
 well-posed in $H^{-3/4}$ with the aid of $I$-method  and the dyadic
 bilinear estimates at the same time. When $\alpha=1$ and $k=2,$ (\ref{1.02})
  reduces to the stochastic
 KdV equation  which has been studied by some people, we refer
 the readers to \cite{Bouard-1998,Bouard-1999,Bouard-2007}.
Recently, motivated by \cite{Bouard-1999}, Chen et al. \cite{CGG}
studied the Cauchy  problem for the stochastic Camassa-Holm equation.

  When $\alpha=0$ and $k=2,$ (\ref{1.03})  reduces to the Benjamin-Ono equation
  which has been studied by many people, we refer the readers to
  \cite{Koch,KT,MM,MA,MS,MF,MP,Tao}. By using the gauge
   transformation introduced by \cite{Tao} and a new bilinear
   estimate, Ionescu and Kenig \cite{IK} proved that the Benjamin-Ono
   equation is globally well-posed in $H^{s}(\R)$  with $s\geq 0.$

When $0<\alpha <1$ and $k=2,$ (\ref{1.03})  has been investigated by some people, we refer the readers to
\cite{CKS,GV,Herr, H,KPV1991}. In \cite{H}, the author proved that
(\ref{1.03})  is locally well-posed in $H^{(s,a)}$ ,$a=\frac{1}{\alpha+1}-\frac{1}{2},s>-\frac{3\alpha}{4}$ and globally well-posed in $H^{(0,a)}$, $a=\frac{1}{\alpha+1}-\frac{1}{2}$.
Recently, by using   a frequency dependent renormalization method, Herr et al.
\cite{HerrCPDE} proved that (\ref{1.03})  is globally well-posed in $L^{2}$ if $0<\alpha <1$ and $k=2$.
Very recently, Guo \cite{GJDE} proved that
(\ref{1.03}) is locally well-posed in $H^{s}$ with $s\geq 1-\alpha$
if $0\leq \alpha \leq 1$ with $k=2$ and in $H^{s}$ with $s\geq \frac{1}{2}-\frac{\alpha}{4},k=3.$

When $\alpha = 1$ and $k = 3$, (\ref{1.03}) reduces to the mKdV equation which has been investigated by
many authors, for instance, see \cite{CKSTJAMS, GIMRN, GJMPA, KPVIUMJ, KPV1993, Kis, TT, N, NTT, KPV2001} and the references therein.
In \cite{Koch}, by using the inverse scattering method, Koch and Tzvetkov proved that the Cauchy
problem for the mKdV equation is locally well-posed on $\mathbf {T}$ in $H^{s}$ with $s \geq0.$ In \cite{TT}, Takaoka
and Tsutsumi proved that the Cauchy problem for the mKdV possesses a unique solution on T
in $H^{s}$ with $\frac{3}{8}<s<\frac{1}{2}$. By using the modified Fourier restriction norm method, Nakanishi et al.\cite{NTT} proved
that the Cauchy problem for the mKdV on $\mathbf{T}$ in $H^{s}$ with $s>\frac{1}{3}$ is locally well-posed and is locally
well-posed in $H^{s}$ with $s>\frac{1}{4}$ with the help of the additional assumption on initial data. Recently,
Molinet \cite{Molinet}  proved that the solution-maps associated with the mKdV equation is discontinuous for
the $H^{s}$   topology for $s < 0$. Soonsik and Oh \cite{SO} studied the unconditional well-posedness of mKV
equation. By using the It$\hat{o}$ formula, BDG inequality and the conserved laws of the KdV equation,
de Bouard and Debussche \cite{Bouard-1998} studied the existence of and uniqueness of solutions to the Cauchy
problem for the Stochastic KdV in $H^{1}(R)$ in the case of additive noise and existence of martingale
solutions in $L^{2}(R)$ in the case of multiplicative noise with the aid of Strichartz estimates and It$\hat{o}$
formula as well as BDG inequality. de Bouard et al. \cite{Bouard-1999} obtained the existence of the solution
to the stochastic KdV in $L^{2}$ with the aid of the modified Bourgain spaces.

In this paper, inspired by \cite{Bouard-1998, Bouard-1999},  we focus on the  case $0<\alpha \leq 1$ and $k=3$ of
 (\ref{1.01}).  By using the  Sobolev spaces and the
Bourgain spaces, we proved that (\ref{1.01})  is locally
well-posed for the initial data $u_{0}(x,w)\in L^{2}
(\Omega; H^{s}(\R))$  with $s\geq\frac{1}{2}-\frac{\alpha}{4}$, where $0< \alpha \leq 1.$
In particular, when $\alpha=1,$ we prove that it is globally well-posed for the initial data
 $u_{0}(x,w)\in L^{2}
(\Omega; H^{1}(\R))$.
Compared to the deterministic KdV and Benjamin-Ono equation, the structure of stochastic Benjamin-Ono equation is more complicated.  The perturbation of the noise destroyed the structure of original structure of Benjamin-Ono.
More precisely, Lemma 2.6 requires $0<b<\frac{1}{2}.$ By using the idea of \cite{Ta}, we firstly  establish the bilinear estimate, then, apply the bilinear estimate which is just Theorem 3.1 to  establish the  trilinear estimate which are Lemmas 4.1-4.2,  thus,  we need to use Lemmas 2.2, 2.3 which are not used in the deterministic KdV and Benjamin-Ono  to establish bilinear and trilinear estimates. Then, the trilinear estimate in combination with the fixed point argument yields Theorem 1.1. For the Theorem 1.2, we use the frequency truncated technique rather than the method of \cite{Bouard-1999}.

We give some notations before giving the main result.
 We denote $X\sim Y$ by
$A_{1}|X|\leq |Y|\leq A_{2}|X| ,$
 where $A_{j}>0\>(j=1,2)$ and denote $X\gg Y$ by $|X|>C|Y|,$
 where $C$ is some positive number which is larger than 2.
$\langle\xi\rangle ^{s} =(1+\xi^{2})^{\frac{s}{2}}$ for any $\xi\in
\R,$ and $\mathscr{F}u$ denotes the Fourier  transformation of $u$
with respect to its all variables. $\mathscr{F}^{-1}u$ denotes the
Fourier inverse transformation of $u$ with respect to its all
variables. $\mathscr{F}_{x}u$ denotes the Fourier  transformation of
$u$ with respect to its  space variable.
$\mathscr{F}^{-1}_{x}u$ denotes the Fourier inverse transformation
of $u$ with respect to its  space variable.
 $H^{s}(\R)$ is the Sobolev space with norm
 $\|f \|_{H^{s}(\SR)}=
 \|\langle\xi\rangle ^{s}\mathscr{F}_{x}{f}\|_{L_{\xi}^{2}(\SR)}$.
For any $s,b\in \R,\> X _{s,\> b}(\R^{2})$  is the Bourgain
space with phase function $\phi(\xi)=\xi|\xi|^{1+\alpha}$. That is, a
function $u(x,t)$    belongs to
$ X _{s, b}(\R^{2})$ iff
$$
 \|u\|_{X _{s, \>b}(\SR^{2})}=
 \left \| \langle\xi\rangle ^{s}\langle \tau-\xi|\xi|^{\alpha+1}\rangle^{b}
  \mathscr{F}{u}(\xi,\tau)\right\| _{{L_{\tau}^{2}(\SR)}{L_{\xi}^{2}(\SR)}}<\infty.
$$
 For any given interval $L$,  $ X_{s,\>
b}(\R\times L)$ is the space of the restriction of all
functions in $ X_{s,\>b}(\R^{2})$ on $\R \times L$, and for $u\in
X_{s, \>b}(\R\times L)$ its norm is
$$
 \| u \|_{ X_{s,\>b}(\SR\times L)}
 =\inf\{\| U\| _{X_{s, \>b}(\SR^{2})};U|_{\SR\times L}=u\}.
$$
When $L=[0, T]$,  $X_{s,\>b}(\R\times L)$ is abbreviated as
$X_{s, b}^{T}$. Throughout this paper, we always assume that
$w(\xi)=\xi|\xi|^{\alpha+1},$
$\psi$ is a
smooth function, $\psi_{\delta}(t)=\psi(\frac{t}{\delta}),$
satisfying $0\leq\psi \leq 1,$  $\psi=1$ when $t\in [0,1],$ ${\rm
supp}\psi\subset[-1,2]$ and $\sigma=\tau-\xi|\xi|^{\alpha+1},$
$\sigma_{k}=\tau_{k}-\xi_{k}|\xi_{k}|^{\alpha+1}\> (k=1,2),$
\begin{eqnarray*}
U(t)u_{0}&=&\frac{1}{\sqrt{2\pi}}\int_{\SR} e^{i(x\xi-t\xi|\xi|^{\alpha+1})}\mathscr{F}_{x}{u}_{0}(\xi)d\xi,\\
\|f\|_{L_{t}^{q}L_{x}^{p}}&=&\left(\int_{\SR}
\left(\int_{\SR}|f(x,t)|^{p}dx\right)^{\frac{q}{p}}dt\right)^{\frac{1}{q}},\\
\|f\|_{L_{t}^{p}L_{x}^{p}}&=&\|f\|_{L_{xt}^{p}}.
\end{eqnarray*}
We assume that $B(x,t)$, $t\geq 0, x\in \R$, is a zero  mean
gaussian process  whose covariance function is given by
\begin{eqnarray*}
{\bf E}(B(t,x)B(s,y))=(t\wedge s)(x\wedge y)
\end{eqnarray*}
for $t,s\geq 0, x,y \in \R$.
$(.,.)$ denotes the $L^{2}$ space duality product, i.e.,
$(f,g)=\int_{\SR}f(x)g(x)dx.$
$(\Omega, \mathscr{F},{\bf P})$  is a probability space  endowed with a
filtration  $(\mathscr{F}_{t})_{t\geq 0}$.
 ${\bf E}f=\int_{\Omega }fd{\bf P}.$
  $W(t)$
 is a cylindrical Wiener process $(W(t))_{t\geq 0}$ on $L^{2}(\R)$
 associated with the filtration $(\mathscr{F_{t}})_{t\geq 0}$.
For any  orthonormal  basis $(e_{k})_{k\in {\bf N}}$ of $L^{2}(\R)$,
$W=\sum_{k=0}^{\infty}\beta_{k}e_{k}$
for a sequence $(\beta_{k})_{k\in {\bf N}}$ of real, mutually
 independent brownian motions on $(\Omega, \mathscr{F},{\bf P},\mathscr{F}_{t})_{t\geq 0}).$
Let $H$ be a Hilbert space, $L_{2}^{0}(L^{2}(\R), H)$ the space of
 Hilbert-Schmidt operators from $L^{2}(\R)$  into $H$. Its norm is given by
$
\|\Phi\|_{L_{2}^{0}(L^{2}(\SR),H)}^{2}=\sum\limits_{j\in {\bf N}}|\Phi e_{j}|_{H}^{2}.
$
When $H=H^{s}(\R)$, $L_{2}^{0}(L^{2}(\R),H^{s}(\R))=L_{2}^{0,s}$.

The main results of this paper  are as follows:

\begin{Theorem}\label{Thm1}
Let   $u_{0}(x,\omega)\in  L^{2}
(\Omega; H^{s}(\R)) $ with $s\geq \frac{1}{2}-\frac{\alpha}{4}$ and $\Phi \in L_{2}^{0,\>s}$ and $u_{0}$ be $\mathscr{F}_{0}$  measurable. Then,
for a.e. $\omega\in \Omega$, there exists a $T_{\omega}>0$ and
 a unique solution of the Cauchy problem for
  (\ref{1.01}) on $[0, T_{\omega}]$ satisfying
  \begin{eqnarray*}
 u\in C([0, T_{\omega}];H^{s}(\R)))\cap X_{s,b}^{T_{\omega}}.
  \end{eqnarray*}
\end{Theorem}
\begin{Theorem}\label{Thm2}
Let  $\alpha=1$, $u_{0}(x,\omega)\in  L^{2}
(\Omega; H^{1}(\R)) $ and $\Phi \in L_{2}^{0,\>1}$ and $u_{0}$
 and $\mathscr{F}_{0}$ be  measurable. Then the solution to the
  Cauchy problem  for (\ref{1.01}) global and belongs to
\begin{eqnarray*}
L^{2}(\Omega; C([0,T_{0}];H^{1}(\R))
\end{eqnarray*}
for any $T_{0}>0.$

\end{Theorem}
The rest of the paper is organized as follows. In Section 2,
some key interpolation inequalities and preliminary estimates
 are established. In the Section 3, we  establish  bilinear estimate
  with the aid of  Fourier restriction norm method. In Section 4,  we will show the trilinear estimate.
In section 5, we prove Theorem 1.1.  In section 6, we  prove Theorem 1.2.

\bigskip
\bigskip
 \noindent{\large\bf 2. Preliminaries }

\setcounter{equation}{0}

\setcounter{Theorem}{0}

\setcounter{Lemma}{0}
In this section, we give some preliminaries which plays the crucial role in establishing the main theorems.

\setcounter{section}{2}
\begin{Lemma}
Let $\theta\in[0,1]$, $\gamma>0$ and $U_\gamma(t)u_0(x)
=\int_Re^{i(t\phi(\xi)+x\xi)}|\phi''(\xi)|^
{\frac{\gamma}{2}}\mathcal{F}_xu_0(\xi)d\xi$.
Then $$
\|U_{\frac{\theta}{2}}(t)u_0\|_{L_t^qL_x^p}\leq C\|u_0\|_{L_x^2},
$$
where $(p,q)=(\frac{2}{1-\theta},\frac{4}{\theta})$.

\end{Lemma}
For the proof of Lemma 2.1, we refer the
readers to Theorem 2.1 of \cite{KPVIUMJ}.
\begin{Lemma}
Let $b=\frac{1}{2}+\epsilon$, $0<\epsilon\ll 1,$ then
\begin{eqnarray}
\|u\|_{L^4_{xt}}\leq C\|u\|_{X_{0
,\frac{\alpha+3}{2(\alpha+2)}(\frac{1}{2}+\epsilon)}}\label{2.01}
\end{eqnarray}
and
\begin{eqnarray}
\left\|D_{x}^{\frac{\alpha}{8}}u\right\|_{L_{xt}^{6}}
\leq C\|u\|_{X_{0,\frac{3}{4}b}}\label{2.02}.
\end{eqnarray}
\end{Lemma}

\begin{proof}
Let $\theta=\frac{2}{3}$, it follows from Lemma 2.1 that
$$
\left\|\int_R e^{it\phi(\xi)+ix\xi}|\phi''(\xi)|^{\frac{1}{6}}
\mathcal{F}_x u_0(\xi)d\xi\right\|_{L^6_{xt}}\leq C\|u_0\|_{L_{\xi}^2}.
$$
where $|\phi|=|\xi|^{\alpha+1}$, $|\phi''|=c|\xi|^{\alpha}$, then
$$
\left\|\int_R e^{it\phi(\xi)+ix\xi}|\xi|^{\frac{\alpha}{6}}
\mathcal{F}_x u_0(\xi)d\xi\right\|_{L^6_{xt}}\leq C\|u_0\|_{L^2_{\xi}}.
$$
Due to $\|f\|_{L_{xt}^{2\alpha+6}}\leq C\|D_x^{\gamma}
D_t^{\gamma}f\|_{L_{xt}^6}$ where $\gamma=\frac{\alpha}{6(\alpha+3)}$. Then
\begin{eqnarray}
\|U(t)u_0(x)\|_{L_{xt}^{2\alpha+4}}&=& C\left\|\int_Re^{i(t\phi+x\xi)}
\mathcal{F}_xu_0(\xi)d\xi\right\|_{L_{xt}^{2\alpha+4}}\nonumber\\
&\leq &C \left\|D_x^{\gamma}D_t^{\gamma}\int_Re^{i(t\phi+x\xi)}
\mathcal{F}_xu_0(\xi)d\xi\right\|_{L_{xt}^6}\nonumber\\
&=& C\left\|\int_Re^{i(t\phi+x\xi)}|\xi|^{\frac{\alpha}{6}}
\mathcal{F}_xu_0(\xi)d\xi\right\|_{L_{xt}^6}\leq C\|u_0\|_{L_x^2}.\label{2.03}
\end{eqnarray}
 Combining  $\|U(t)u_0(x)\|_{L_{xt}
^{2\alpha+6}}\leq C\|u_0\|_{L_x^2}$  with a standard argument,  we have
\begin{eqnarray}
\|u(x)\|_{L_{xt}^{2\alpha+6}}\leq C\|u\|_{X_0,\frac{1}{2}+\epsilon}.\label{2.04}
\end{eqnarray}
By using the  Plancherel identity,  we have that
\begin{eqnarray}
\|u\|_{L_{xt}^2}=C\|u\|_{X_{0,0}}.\label{2.05}
\end{eqnarray}
Interpolating  (\ref{2.04})   with  (\ref{2.05})  yields
\begin{eqnarray}
\|u\|_{L_{xt}^4}\leq C\|u\|_{X_{0,\frac{\alpha+3}{2(\alpha+2)}
(\frac{1}{2}+\epsilon)}}.\label{2.06}
\end{eqnarray}
From (\ref{2.03}),  by  using  a standard proof,  we have that
\begin{eqnarray}
\|D_{x}^{\frac{\alpha}{6}}u\|_{L_{xt}^{6}}\leq
C\|u\|_{X_{0,b}}.\label{2.07}
\end{eqnarray}
Interpolating (\ref{2.07})  with  (\ref{2.05})  yields
\begin{eqnarray}
\|D_{x}^{\frac{\alpha}{8}}u\|_{L_{xt}^{4}}\leq
C\|u\|_{X_{0,\frac{3}{4}b}}.\label{2.08}
\end{eqnarray}

We have completed the proof of  Lemma 2.2.

\end{proof}

\begin{Lemma}

Let $b=\frac{1}{2}+\epsilon.$ Then,  for $0\leq s\leq \frac{1}{2},$  we have that
\begin{eqnarray}
\left\|I^{s}(u_{1},u_{2})\right\|
  _{L_{xt}^{2}}\leq C\prod_{j=1}^{2}\| u_{j}\| _{X_{0,
  \frac{\alpha+3+2(\alpha+1) s}{2(\alpha+2)}b}},\label{2.09}
\end{eqnarray}
where
\begin{eqnarray*}
\mathscr{F}I^{s}(u_{1},u_{2})(\xi,\tau)
&=&\int_{\!\!\!\mbox{\scriptsize $
\begin{array}{l}
\xi=\xi_{1}+\xi_{2}\\
\tau=\tau_{1}+\tau_{2}
\end{array}
$}}
||\xi_{1}|^{\alpha+1}-|\xi_{2}|^{\alpha+1}|^{s}
\mathscr{F}{u_{1}}(\xi_{1},\tau_{1})\mathscr{F}{u_{2}}
(\xi_{2},\tau_{2})\,d\xi_{1}d\tau_{1}.
\end{eqnarray*}
\end{Lemma}
{\bf Proof.}
  Let $F_{j}(\xi_{j},\tau_{j})=\langle\sigma_{j} \rangle^{
\frac{\alpha+3+2(\alpha+1) s}{2\alpha+4}b}\mathscr{F}{u_{j}}(\xi_{j},\tau_{j})(j=1,2)$.
To prove Lemma 2.3, by the  Plancherel identity,  it suffices to prove that
\begin{eqnarray}
&&\left\|\int_{\!\!\!\mbox{\scriptsize $
\begin{array}{l}
\xi=\xi_{1}+\xi_{2}\\
\tau=\tau_{1}+\tau_{2}
\end{array}
$}}||\xi_{1}|^{\alpha+1}-|\xi_{2}|^{\alpha+1}|^{s}\frac{F_{1}}{\langle\sigma_{1} \rangle^{
\frac{\alpha+3+2(\alpha+1) s}{2\alpha+2}b}}\frac{F_{2}}{\langle\sigma_{2} \rangle^{
\frac{\alpha+3+2(\alpha+1) s}{2\alpha+2}b}}d\xi_{1}d\tau_{1}\right\|_{L_{\xi\tau}^{2}}\nonumber\\&&\leq C
\prod_{j=1}^{2}\|F_{j}\|_{L_{\xi\tau}^{2}}.\label{2.010}
\end{eqnarray}
Assume that $b_{1}=\frac{\alpha+3+2(\alpha+1)s}{2\alpha+4}b$. By using the Young inequality,
since $0<s<\frac{1}{2},$ we have that
\begin{eqnarray}
&&||\xi_{1}|^{\alpha+1}-|\xi_{2}|^{\alpha+1}|^{s}\langle\sigma_{1}
\rangle^{ -b_{1}}\langle\sigma_{2}
\rangle^{ -b_{1}}\nonumber\\&&=||\xi_{1}|^{\alpha+1}-|\xi_{2}|^{\alpha+1}|^{s}
\langle\sigma_{1} \rangle^{ -2bs}\langle\sigma_{2}
\rangle^{ -2bs}\langle\sigma_{1} \rangle^{ -(b_{1}-2bs)}\langle\sigma_{2}
\rangle^{ -(b_{1}-2bs)}\nonumber\\&&\leq2s||\xi_{1}|^{\alpha+1}-|\xi_{2}|^{\alpha+1}|^{1/2}
\langle\sigma_{1} \rangle^{ -b}\langle\sigma_{2} \rangle^{ -b}+(1-2s)
\langle\sigma_{1} \rangle^{ -\frac{\alpha+3}{2\alpha+4}b}\langle\sigma_{2} \rangle^{
-\frac{\alpha+3}{2\alpha+4}b}
\nonumber\\&&\leq||\xi_{1}|^{\alpha+1}-|\xi_{2}|^{\alpha+1}|^{1/2}\langle\sigma_{1}
\rangle^{ -b}\langle\sigma_{2} \rangle^{ -b}+\langle\sigma_{1} \rangle^{
 -\frac{\alpha+3}{2\alpha+4}b}\langle\sigma_{2}
 \rangle^{ -\frac{\alpha+3}{2\alpha+4}b}\label{2.011}.
\end{eqnarray}
By using (\ref{2.011}),  Plancherel identity,
 Lemma 3.1 in \cite{H}, we have that
\begin{eqnarray}
&&\left\|\int_{\!\!\!\mbox{\scriptsize $
\begin{array}{l}
\xi=\xi_{1}+\xi_{2}\\
\tau=\tau_{1}+\tau_{2}
\end{array}
$}}||\xi_{1}|^{\alpha+1}-|\xi_{2}|^{\alpha+1}|^{s}\frac{F_{1}}{\langle\sigma_{1} \rangle^{
\frac{\alpha+3+2(\alpha+1)s}{2\alpha+4}b}}\frac{F_{2}}{\langle\sigma_{2} \rangle^{ \frac{\alpha+3+2(\alpha+1) s}{2\alpha+2}b}}
d\xi_{1}d\tau_{1}\right\|_{L_{\xi\tau}^{2}}\nonumber\\&&\leq
\left\|\int_{\!\!\!\mbox{\scriptsize $
\begin{array}{l}
\xi=\xi_{1}+\xi_{2}\\
\tau=\tau_{1}+\tau_{2}
\end{array}
$}}||\xi_{1}|^{\alpha+1}-|\xi_{2}|^{\alpha+1}|^{1/2}\prod_{j=1}^{2}\frac{F_{j}}{\langle\sigma_{j}
\rangle^{ b}}d\xi_{1}d\tau_{1}\right\|_{L_{\xi\tau}^{2}}\nonumber\\&&\qquad
+\left\|\int_{\!\!\!\mbox{\scriptsize $
\begin{array}{l}
\xi=\xi_{1}+\xi_{2}\\
\tau=\tau_{1}+\tau_{2}
\end{array}
$}}\prod_{j=1}^{2}\frac{F_{1}}{\langle\sigma_{j} \rangle^{ \frac{\alpha+3}{2\alpha+4}b}}
d\xi_{1}d\tau_{1}\right\|_{L_{\xi\tau}^{2}}\nonumber\\
&&\leq C\prod_{j=1}^{2}\left\|\mathscr{F}^{-1}\left(\frac{F_{j}}{\langle\sigma_{j} \rangle^{ b}}\right)\right\|_{X_{0,b}}+C\prod_{j=1}^{2}\left\|\mathscr{F}^{-1}\left(\frac{F_{j}}{\langle\sigma_{j}
 \rangle^{ \frac{\alpha+3}{2\alpha+4}b}}\right)\right\|_{X_{0,\frac{\alpha+3}{2\alpha+4}b}}
\nonumber\\
&&\leq C\prod_{j=1}^{2}\|F_{j}\|_{L_{\xi\tau}^{2}}.\label{2.012}
\end{eqnarray}

We have completed the proof of Lemma 2.3.

\begin{Lemma} \label{Lemma 2.4}
Let $u_{0}\in H^{s}(\R)$, $c>1/2,0<b<1/2.$ Then for $t\in [0,T]$, $U(t)u_{0}\in X_{s,\>c}^{T}$ and there is a constant $k_{2}>0$ such that
\begin{eqnarray}
&&\|U(t)u_{0}\|_{X_{s,c}^{T}}\leq k_{2}\|u_{0}\|_{H^{s}}.\label{2.013}
\end{eqnarray}
There is a constant $c>0$ such that for $t\in [0,1]$ and $f\in X_{s,\>b}^{T},$
\begin{eqnarray}
\left\|\int_{0}^{T}U(t-s)f(s)ds\right\|_{X_{s,b}^{T}}\leq CT^{1-2b}\|f\|_{X_{s,-b}^{T}}.\label{2.014}
\end{eqnarray}
\end{Lemma}

For the proof of Lemma 2.4, we refer the readers to  Lemma 3.1 of \cite{Bouard-1999}.
\begin{Lemma} \label{Lemma 2.5}
Let
\begin{eqnarray*}
\overline{u}=\int_{0}^{t}U(t-s)\Phi dW(s)
\end{eqnarray*}
and $\Phi \in L_{2}^{0,s}$, for $t\in [0,T],$  we have
\begin{eqnarray}
E(\sup\limits_{t\in[0,T]}\|\overline{u}\|_{H^{s}}^{2})\leq38 T\|\Phi\|_{L_{2}^{0,s}}^{2}.\label{2.015}
\end{eqnarray}
\end{Lemma}

 Lemma 2.5 can be proved similarly to Proposition 2.1 of \cite{Bouard-1999}.

\begin{Lemma}\label{Lemma2.6}
Let
\begin{eqnarray*}
\bar{u}=\int_{0}^{t}U(t-s)\Phi dW(s),
\end{eqnarray*}
$s,b\in \R$ with $b<\frac{1}{2}$ and $\Phi\in  L_{2}^{0,s}$.  Then, we have that
\begin{eqnarray}
{\rm E}\left(\|\psi \bar{u}\|_{X_{s,b}}^{2}\right)\leq C\|\Phi\|_{L_{2}^{0,s}}^{2}.\label{2.018}
\end{eqnarray}
\end{Lemma}

For the proof of Lemma 2.6, we refer the readers to Proposition 2.1 of \cite{Bouard-1999}.
\bigskip
\bigskip

\section{Bilinear estimate}
\setcounter{equation}{0}

 \setcounter{Theorem}{0}

\setcounter{Lemma}{0}

 \setcounter{section}{3}
In this section, we give an important bilinear estimate which can be used to establish two important trilinear estimates.
\begin{Theorem}
For all $u,v$ on $\R\times \R$, $0\ll \epsilon\leq 1$ and $b=\frac{1}{2}-\epsilon$, we have
\begin{eqnarray}
\|u_{1}u_{2}\|_{L^{2}}\leq C\|u_{1}\|_{X_{-\frac{1}{2},b}}\|u_{2}\|_{X_{\frac{1}{2}-\frac{\alpha}{4},b}}.\label{3.01}
\end{eqnarray}
\end{Theorem}
\noindent {\bf Proof.}
Define
\begin{eqnarray*}
&&F_1(\xi_1,\tau_1)=\langle \xi_{1}\rangle^{-1/2}\langle\sigma_1\rangle^{b}\mathcal{F}u_1(\xi_1,\tau_1)
\quad F_{2}(\xi_{2},\tau_{2})=\langle\xi_{2}\rangle^{\frac{1}{2}-\frac{\alpha}{4}}\langle\sigma_{2}\rangle^{b}
\mathcal{F}u(\xi_{2},\tau_{2})\nonumber\\
&&\sigma_j=\tau_j-|\xi_{j}|^{\alpha+1}\xi_{j},\;j=1,2.
\end{eqnarray*}
To obtain (\ref{3.01}), it suffices to prove that
\begin{eqnarray}
\int_{\SR^2}\int_{\tiny\begin{array}{cc}\xi=\xi_1+\xi_2\\ \tau =\tau_1+\tau_2\end{array}}K_{1}(\xi_1,\tau_1,\xi,\tau)|F|\prod_{j=1}^2|F_j|
d\xi_1d\tau_1d\xi d\tau
\leq C\|F\|_{L_{\xi\tau}^2}\prod_{j=1}^{2}\|F_j\|_{L_{\xi\tau}^2},\label{3.2}
\end{eqnarray}
where
\begin{eqnarray*}
K_{1}(\xi_1,\tau_1,\xi,\tau)=\frac{\langle \xi_{1}\rangle^{1/2}
\langle\xi_{2}\rangle^{\frac{\alpha}{4}-\frac{1}{2}}}{\langle\sigma_{1}\rangle^{b}\langle\sigma_{2}\rangle^{b}}.
\end{eqnarray*}
Without loss of  generality,  we assume that  $F\geq 0,F_j\geq 0(j=1,2).$

\begin{eqnarray*}
&&\hspace{-0.8cm}\Omega_1=\{(\xi_1,\tau_1,\xi,\tau)\in {\rm R^4},\xi=\sum_{j=1}^{2}\xi_j,\tau=\sum_{j=1}^{2}\tau_j,
|\xi_1|\leq |\xi_2|\leq 6 \},\\
&&\hspace{-0.8cm} \Omega_2=\{ (\xi_1,\tau_1,\xi,\tau)\in {\rm R^4},\xi=\sum_{j=1}^{2}\xi_j,\tau=\sum_{j=1}^{2}\tau_j,
|\xi_2|\geq 6, |\xi_{2}|\gg|\xi_{1}|\},\\
&&\hspace{-0.8cm}\Omega_3=\{(\xi_1,\tau_1,\xi,\tau)\in {\rm R^4},\xi=\sum_{j=1}^{2}\xi_j,\tau=\sum_{j=1}^{2}\tau_j,|\xi_2|\geq6,
|\xi_2|\sim |\xi_1|\},\\
&&\hspace{-0.8cm}\Omega_4=\{(\xi_1,\tau_1,\xi,\tau)\in {\rm R^4},\xi=\sum_{j=1}^{2}\xi_j,\tau=\sum_{j=1}^{2}\tau_j,
|\xi_2|\leq |\xi_1|\leq 6 \},\\
&&\hspace{-0.8cm}\Omega_5=\{(\xi_1,\tau_1,\xi,\tau)\in {\rm R^4},\xi=\sum_{j=1}^{2}\xi_j,\tau=\sum_{j=1}^{2}\tau_j,
|\xi_1|\geq6,|\xi_{1}|\gg|\xi_{2}|\},\\
&&\hspace{-0.8cm}\Omega_6=\{(\xi_1,\tau_1,\xi,\tau)\in {\rm R^4},\xi=\sum_{j=1}^{2}\xi_j,\tau=\sum_{j=1}^{2}\tau_j,|\xi_1|\geq6,|\xi_{1}|\geq |\xi_{2}|,
|\xi_1|\sim |\xi_2|\},\\
\end{eqnarray*}
  We define
\begin{eqnarray*}
f_{j}=\mathscr{F}^{-1}
\frac{F_{j}}{\langle\sigma_{j}\rangle^{b}},j=1,2.
\end{eqnarray*}
\noindent
(1).
$
\Omega_1=\{(\xi_1,\tau_1,\xi,\tau)\in {\rm R^4},\xi=\sum_{j=1}^{2}\xi_j,\tau=\sum_{j=1}^{2}\tau_j,
|\xi_1|\leq |\xi_2|\leq 6 \}.
$
In this subregion, we have that
\begin{eqnarray*}
K_{1}(\xi_{1},\tau_{1},\xi,\tau)\leq \frac{C}{\prod_{j=1}^{2}\langle\sigma\rangle^{b}}.
\end{eqnarray*}
By using the Plancherel identity and  the H\"older inequality and  $\frac{\alpha+3}{2(\alpha+2)}(\frac{1}{2}+\epsilon)<\frac{1}{2}-\epsilon$,
we have that
\begin{eqnarray}
&&J_{1}\leq C\int_{\SR^2}\int_{\tiny\begin{array}{cc}\xi=\xi_1+\xi_2\\ \tau =\tau_1+\tau_2\end{array}}\frac{F\prod_{j=1}^{2}F_{j}}{\prod_{j=1}^{2}\langle\sigma_{j}\rangle^{b}}d\xi_{1}d\tau_{1}d\xi d\tau\nonumber\\
&&\leq C\int_{\SR^2}\mathscr{F}^{-1}(F)f_{1}f_{2}dxdt\leq C\|\mathscr{F}^{-1}(F)\|_{L_{xt}^{2}}\prod_{j=1}^{2}\|f_{j}\|_{L_{xt}^{4}}\nonumber\\
&&\leq C\|F\|_{L_{\xi\tau}^{2}}\prod_{j=1}^{2}\|f_{j}\|_{X_{0,\frac{\alpha+3}{2(\alpha+2)}(\frac{1}{2}+\epsilon)}}\nonumber\\
&&\leq C\|F\|_{L_{\xi\tau}^{2}}\prod_{j=1}^{2}\|F_{j}\|_{L_{\xi\tau}^{2}}.
\end{eqnarray}

\noindent(2).$\Omega_2=\{ (\xi_1,\tau_1,\xi,\tau)\in {\rm R^4},\xi=\sum_{j=1}^{2}\xi_j,\tau=\sum_{j=1}^{2}\tau_j,
|\xi_2|\geq 6, |\xi_{2}|\gg|\xi_{1}|\}.$

\noindent
If $|\xi_{1}|\leq 1,$ we have that
\begin{eqnarray*}
K_{1}(\xi_{1},\tau_{1},\xi,\tau)\leq \frac{C}{\prod_{j=1}^{2}\langle\sigma_{j}\rangle^{b}}
\end{eqnarray*}
This case can be proved similarly to $\Omega_1.$

\noindent If $|\xi_{1}|\geq 1,$ we have
\begin{eqnarray*}
K_{1}(\xi_{1},\tau_{1},\xi,\tau)\leq C\frac{|\xi_{2}|^{\frac{\alpha}{4}}}{\prod_{j=1}^{2}\langle\sigma_{j}\rangle^{b}}
\leq C\frac{||\xi_{2}|^{\alpha+1}-|\xi_{2}|^{\alpha+1}|^{\frac{\alpha}{4(\alpha+1)}}}{\prod_{j=1}^{2}\langle\sigma_{j}\rangle^{b}}.
\end{eqnarray*}
By using Lemma 2.3, we have
\begin{eqnarray*}
&&J_{2}\leq C\int_{\SR^2}\int_{\tiny\begin{array}{cc}\xi=\xi_1+\xi_2\\ \tau =\tau_1+\tau_2\end{array}}\frac{||\xi_{2}|^{\alpha+1}-|\xi_{2}|^{\alpha+1}|^{\frac{\alpha}{4(\alpha+1)}}
F\prod_{j=1}^{2}F_{j}}{\prod_{j=1}^{2}\langle\sigma_{j}\rangle^{b}}d\xi_{1}d\tau_{1}d\xi d\tau\nonumber\\
&&\leq C \|F\|_{L_{\xi\tau}^{2}}\left\|\int_{\tiny\begin{array}{cc}\xi=\xi_1+\xi_2\\ \tau =\tau_1+\tau_2\end{array}}\frac{||\xi_{2}|^{\alpha+1}-|\xi_{2}|^{\alpha+1}|^{\frac{\alpha}{4(\alpha+1)}}F\prod_{j=1}^{2}F_{j}}{\prod_{j=1}^{2}
\langle\sigma_{j}\rangle^{b}}d\xi_{1}d\tau_{1}\right\|_{L_{\xi\tau}^{2}}\nonumber\\
&&\leq C\|F\|_{L_{\xi\tau}^{2}}\prod_{j=1}^{2}\|F_{j}\|_{L_{\xi\tau}^{2}}.
\end{eqnarray*}

\noindent(3).
$
\Omega_3=\{(\xi_1,\tau_1,\xi,\tau)\in {\rm R^4},\xi=\sum_{j=1}^{2}\xi_j,\tau=\sum_{j=1}^{2}\tau_j,|\xi_2|\geq6,
|\xi_2|\sim |\xi_1|\}.$
\begin{eqnarray*}
K_{1}(\xi_{1},\tau_{1},\xi,\tau)\leq C\frac{|\xi_{2}|^{\frac{\alpha}{4}}}{\prod_{j=1}^{2}\langle\sigma_{j}\rangle^{b}}\leq C\frac{\prod_{j=1}^{2}|\xi_{j}|^{\frac{\alpha}{8}}}{\prod_{j=1}^{2}\langle\sigma_{j}\rangle^{b}}
\end{eqnarray*}
By using the Plancherel identity and the Cauchy-Schwartz inequality, we have
\begin{eqnarray}
&&J_{3}\leq C\int_{\SR^2}\int_{\tiny\begin{array}{cc}\xi=\xi_1+\xi_2\\ \tau =\tau_1+\tau_2\end{array}}\frac{\prod_{j=1}^{2}|\xi_{j}|^{\frac{\alpha}{8}}
F\prod_{j=1}^{2}F_{j}}{\prod_{j=1}^{2}\langle\sigma_{j}\rangle^{b}}d\xi_{1}d\tau_{1}d\xi d\tau\nonumber\\
&&\leq C \|F\|_{L_{\xi\tau}^{2}}\prod_{j=1}^{2}\left\|D_{x}^{\frac{\alpha}{8}}
\mathscr{F}^{-1}\left(\frac{F_{j}}{\langle\sigma_{j}\rangle^{b}}\right) \right\|_{L_{xt}^{4}}\nonumber\\
&&\leq C\|F\|_{L_{\xi\tau}^{2}}\prod_{j=1}^{2}\|F_{j}\|_{L_{\xi\tau}^{2}}.
\end{eqnarray}

\noindent(4).
$\Omega_4=\{(\xi_1,\tau_1,\xi,\tau)\in {\rm R^4},\xi=\sum_{j=1}^{2}\xi_j,\tau=\sum_{j=1}^{2}\tau_j,
|\xi_2|\leq |\xi_1|\leq 6 \}.$
In this subregion, we have that
\begin{eqnarray*}
K_{1}(\xi_{1},\tau_{1},\xi,\tau)\leq \frac{C}{\prod_{j=1}^{2}\langle\sigma_{j}\rangle^{b}}.
\end{eqnarray*}
Thus subregion can be proved similarly to $\Omega_1.$

\noindent(5).
$\Omega_5=\{(\xi_1,\tau_1,\xi,\tau)\in {\rm R^4},\xi=\sum_{j=1}^{2}\xi_j,\tau=\sum_{j=1}^{2}\tau_j,
|\xi_1|\geq6,|\xi_{1}|\gg|\xi_{2}|\}.$ In this subregion, we have
\begin{eqnarray*}
K_{1}(\xi_{1},\tau_{1},\xi,\tau)\leq C\frac{||\xi_{1}|^{\alpha+1}-|\xi_{2}|^{\alpha+1}|^{\frac{1}{2(\alpha+1)}}}{\prod_{j=1}^{2}\langle\sigma_{j}\rangle^{b}}.
\end{eqnarray*}
By using Lemma 2.3, we have that
\begin{eqnarray*}
&&J_{5}\leq C\int_{\SR^2}\int_{\tiny\begin{array}{cc}\xi=\xi_1+\xi_2\\ \tau =\tau_1+\tau_2\end{array}}\frac{||\xi_{2}|^{\alpha+1}-|\xi_{1}|^{\alpha+1}|^{\frac{1}{2(\alpha+1)}}
F\prod_{j=1}^{2}F_{j}}{\prod_{j=1}^{2}\langle\sigma_{j}\rangle^{b}}d\xi_{1}d\tau_{1}d\xi d\tau\nonumber\\
&&\leq C \|F\|_{L_{\xi\tau}^{2}}\left\|\int_{\tiny\begin{array}{cc}\xi=\xi_1+\xi_2\\ \tau =\tau_1+\tau_2\end{array}}\frac{||\xi_{2}|^{\alpha+1}-|\xi_{1}|^{\alpha+1}|^{\frac{1}{2(\alpha+1)}}F\prod_{j=1}^{2}F_{j}}{\prod_{j=1}^{2}
\langle\sigma_{j}\rangle^{b}}d\xi_{1}d\tau_{1}\right\|_{L_{\xi\tau}^{2}}\nonumber\\
&&\leq C\|F\|_{L_{\xi\tau}^{2}}\prod_{j=1}^{2}\|F_{j}\|_{L_{\xi\tau}^{2}}.
\end{eqnarray*}

\noindent(6).
$\Omega_6=\{(\xi_1,\tau_1,\xi,\tau)\in {\rm R^4},\xi=\sum_{j=1}^{2}\xi_j,\tau=\sum_{j=1}^{2}\tau_j,|\xi_1|\geq6,|\xi_{1}|\geq |\xi_{2}|,
|\xi_1|\sim |\xi_2|\}.$

This subregion can be proved similarly to $\Omega_3.$

We have completed the proof of Lemma 3.1.
\section{Trilinear estimates}
\setcounter{equation}{0}

 \setcounter{Theorem}{0}

\setcounter{Lemma}{0}

 \setcounter{section}{4}
In this section, we will establish two new trilinear estimates which play
a crucial role in establishing the local well-posedness of solution.

We will establish the Lemma 4.1 with the aid of  the idea in \cite{Tao}. Let $Z=\R$  and  $\Gamma_{k}(Z)$
denote the hyperplane in $\R^{k}$
\begin{eqnarray*}
\Gamma_{k}(Z):=\left\{(\xi_{1},\cdot\cdot\cdot,\xi_{k})\in Z^{k}, \xi_{1}+\cdot\cdot\cdot+\xi_{k}=0\right\}
\end{eqnarray*}
endowed with the induced measure
\begin{eqnarray*}
\int_{\Gamma_{k}(Z)}f:=\int_{Z^{k-1}}f(\xi_{1},\cdot\cdot\cdot,\xi_{k-1},-\xi_{1}-\cdot\cdot\cdot-\xi_{k-1})d\xi_{1}\cdot\cdot\cdot d\xi_{k}.
\end{eqnarray*}
A function $m:\Gamma_{k}(Z)\rightarrow C$ is said to be a $[k;Z]$-multiplier, and we define the norm $\|m\|_{[k;Z]}$ to be the best constant such that the inequality
\begin{eqnarray*}
\left|\int_{\Gamma_{k}(Z)}m(\xi)\prod_{j=1}^{k}f_{j}(\xi_{j})\right|\leq \|m\|_{[k;Z]}\prod_{j=1}^{k}\|f_{j}\|_{L^{2}}.
\end{eqnarray*}
holds for all test function $f_{j}$  on $Z.$
\begin{Lemma}
Let $s_{0}= \frac{1}{2}-\frac{\alpha}{4},$
$b=\frac{1}{2}-\epsilon$. Then
\begin{eqnarray}
\|\partial_x(u_1u_2u_3)\|_{X_{s_{0},-b}}\leq C\prod_{j=1}^{3}\|u_j\|_{X_{s_{0},b}}.\label{4.01}
\end{eqnarray}
\end{Lemma}
\noindent {\bf Proof.}
By duality, Plancherel identity and the definition, to obtain (\ref{4.01}), it suffices to prove that
\begin{eqnarray}
\left\|\frac{(\sum_{j=1}^{3}\xi_{j})\langle\xi_{4}\rangle^{\frac{1}{2}-\frac{\alpha}{4}}}
{\prod_{j=1}^{4}\langle\tau_{j}-w(\xi_{j})
\rangle^{\frac{1}{2}-\epsilon}\prod_{j=1}^{3}
\langle\xi_{j}\rangle^{\frac{1}{2}-\frac{\alpha}{4}}}\right\|_{[4;\SR\times \SR]}\leq C\label{4.02}.
\end{eqnarray}
By using the symmetry and
\begin{eqnarray*}
\langle\xi_{4}\rangle^{\frac{3}{2}-\frac{\alpha}{4}}\leq C\langle\xi_{4}\rangle^{\frac{1}{2}}\left[\sum_{j=1}^{3}\langle\xi_{j}\rangle^{1-\frac{\alpha}{4}}\right]
\end{eqnarray*}
resulting from
\begin{eqnarray*}
|\xi_{1}+\xi_{2}+\xi_{3}|\leq \langle\xi_{4}\rangle ,
\end{eqnarray*}
to obtain (\ref{4.02}), it suffices to prove
\begin{eqnarray}
\left\|\frac{\langle\xi_{4}\rangle^{1/2}\langle\xi_{2}\rangle^{1/2}}
{\langle\xi_{1}\rangle^{\frac{1}{2}-\frac{\alpha}{4}}\langle\xi_{3}\rangle^{\frac{1}{2}-\frac{\alpha}{4}}
\prod_{j=1}^{4}\langle\tau_{j}-w(\xi_{j})
\rangle^{\frac{1}{2}-\epsilon}}\right\|_{[4;\SR\times\SR]}\leq C\label{4.03}.
\end{eqnarray}
(\ref{4.03}) follows from $TT^{\star}$ identity in Lemma 3.7 of \cite{Tao} and Lemma 3.1.

We have completed the proof of Lemma 4.1.

\begin{Lemma}
Let $s\geq s_{0}= \frac{1}{2}-\frac{\alpha}{4},$
$b=\frac{1}{2}-\epsilon$. Then
\begin{eqnarray}
\label{bilinear estimate}
\|\partial_x(u_1u_2u_3)\|_{X_{s,-b}}\leq C\prod_{j=1}^{3}\|u_j\|_{X_{s,b}}.\label{4.04}
\end{eqnarray}

\end{Lemma}
\noindent {\bf Proof.}
(\ref{4.04})  is equivalent to the following inequality
\begin{eqnarray}
\int_{\SR^2}\int_{\tiny\begin{array}{cc}\xi=\xi_1+\xi_2+\xi_{3}\\ \tau =\tau_1+\tau_2+\tau_{3}\end{array}}\frac{|\xi|\langle\xi\rangle^{s}F\prod_{j=1}^{3}F_{j}}{\langle\sigma\rangle^{b}
\prod_{j=1}^{3}\langle\xi_{j}\rangle^{s}\langle\sigma_{j}\rangle^{b}}d\xi_{1}d\tau_{1}d\xi_{2}d\tau_{2}d\xi d\tau\leq C
\|F\|_{L_{\xi\tau}^{2}}\prod_{j=1}^{3}\|F_{j}\|_{L_{\xi\tau}^{2}}.\label{4.05}
\end{eqnarray}
Since
\begin{eqnarray}
\langle\xi\rangle^{s-s_{0}}\leq C\prod_{j=1}^{3}\langle\xi_{j}\rangle^{s-s_{0}},\label{4.06}
\end{eqnarray}
(\ref{4.05})  is equivalent to the following inequality
\begin{eqnarray}
\int_{\SR^2}\int_{\tiny\begin{array}{cc}\xi=\xi_1+\xi_2+\xi_{3}\\ \tau =\tau_1+\tau_2+\tau_{3}\end{array}}\frac{|\xi|\langle\xi\rangle^{s_{0}}F\prod_{j=1}^{3}F_{j}}{\langle\sigma\rangle^{b}
\prod_{j=1}^{3}\langle\xi_{j}\rangle^{s_{0}}\langle\sigma_{j}\rangle^{b}}d\xi_{1}d\tau_{1}d\xi_{2}d\tau_{2}d\xi d\tau\leq C
\|F\|_{L_{\xi\tau}^{2}}\prod_{j=1}^{3}\|F_{j}\|_{L_{\xi\tau}^{2}},\label{4.07}
\end{eqnarray}
which is just the Lemma 4.1.

We have completed the proof of Lemma 4.2.

\begin{Lemma}
Let $s\geq s_{0}= \frac{1}{2}-\frac{\alpha}{4},$
$b=\frac{1}{2}-\epsilon$. Then
\begin{eqnarray}
\label{bilinear estimate}
\|\partial_x(u_1u_2u_3)\|_{X_{s,-b}^{T}}\leq C\prod_{j=1}^{3}\|u_j\|_{X_{s,b}^{T}}.\label{4.08}
\end{eqnarray}

\end{Lemma}

Combining Lemma 4.2 with a standard proof, we can obtain Lemma 4.3.
\bigskip
\bigskip
\section{Local well-posedness}
\setcounter{equation}{0}

 \setcounter{Theorem}{0}

\setcounter{Lemma}{0}

 \setcounter{section}{5}
In this section, we prove Theorem 1.1.
\noindent Let $z(t)=U(t)u_{0}$ and $\bar{u}=\int_{0}^{t}U(t-s)\Phi dW$.

 The solution to (\ref{1.01}) is equivalent to the following integral
equation
\begin{eqnarray}
u(t)=U(t)u_{0}+\frac{1}{3}\int_{0}^{t}U(t-s)\partial_{x}(u^{3})ds+\int_{0}^{t}U(t-s)\Phi dW\label{5.02}.
\end{eqnarray}
 and $v(t)=u(t)-z(t)-\bar{u}$. Then, we have that
\begin{eqnarray}
v(t)=u(t)-z(t)-\bar{u}=\frac{1}{3}\int_{0}^{t}U(t-s)\partial_{x}(v+z(t)+\bar{u})^{3}ds\label{5.03}.
\end{eqnarray}
We define
\begin{eqnarray}
G(v)=\frac{1}{3}\int_{0}^{t}U(t-s)\partial_{x}(v+z(t)+\bar{u})^{3}ds\label{5.04}.
\end{eqnarray}
By using Lemma 4.4, Lemmas 2.4, 2.5, 2.7,  we have that
\begin{eqnarray}
&&\|G(v)\|_{X_{s,b}^{T}}\leq \left\|\frac{1}{3}\int_{0}^{t}U(t-s)\partial_{x}(v+z(t)+\bar{u})^{3}ds\right\|_{X_{s,b}^{T}}\nonumber\\&&\leq
CT^{1-2b}\left(\|v\|_{X_{s,b}^{T}}^{3}+\|z(t)\|_{X_{s,b}^{T}}^{3}+\|\psi\left(\frac{t}{T}\right)\bar{u}\|_{X_{s,b}}^{3}\right)\nonumber\\
&&\leq CT^{1-2b}\left(\|v\|_{X_{s,b}^{T}}^{3}+\|u_{0}\|_{X_{s,b}^{T}}^{3}+\|\psi\left(\frac{t}{T}\right)\bar{u}\|_{X_{s,b}}^{3}\right),\label{5.05}
\end{eqnarray}
similarly, we have that
\begin{eqnarray}
&&\|G(v_{1})-G(v_{2})\|_{X_{s,b}^{T}}\leq \left\|\frac{1}{3}\int_{0}^{t}U(t-s)\partial_{x}(v+z(t)+\bar{u})^{3}ds\right\|_{X_{s,b}^{T}}\nonumber\\&&\leq
CT^{1-2b}\|v_{1}-v_{2}\|_{X_{s,b}^{T}}
\left(\|v_{1}\|_{X_{s,b}^{T}}^{2}+\|v_{2}\|_{X_{s,b}^{T}}^{2}+\|z(t)\|_{X_{s,b}^{T}}^{2}+\|\psi\left(\frac{t}{T}\right)\bar{u}\|_{X_{s,b}}^{2}\right)\nonumber\\&&\leq C
T^{1-2b}\|v_{1}-v_{2}\|_{X_{s,b}^{T}}
\left(\|v_{1}\|_{X_{s,b}^{T}}^{2}+\|v_{2}\|_{X_{s,b}^{T}}^{2}+\|u_{0}\|_{H^{s}}^{2}+\|\psi\left(\frac{t}{T}\right)\bar{u}\|_{X_{s,b}}^{2}\right),\label{5.06}
\end{eqnarray}
Let
\begin{eqnarray}
R_{\omega}=\left[\|\psi\left(\frac{t}{T}\right)\bar{u}\|_{X_{s,b}}+\|u_{0}\|_{H^{s}}+2\right]^{3}.\label{5.07}
\end{eqnarray}
and  define
\begin{eqnarray}
T_{\omega}={\rm inf}\left\{T>0, CT^{1-2b}R_{\omega}^{3}\geq \frac{1}{4}\right\}.\label{5.08}
\end{eqnarray}
From Lemma 2.6, for any $0<T<1,$  we have that
\begin{eqnarray*}
\|\chi_{[0,T]}\bar{u}\|_{X_{s,b}}\leq C\|\bar{u}\|_{X_{s,b}^{1}}\leq C(\omega)
\end{eqnarray*}
a.s. Moreover, since $b=\frac{1}{2}-\epsilon$, $\|\chi_{[0,T]}\bar{u}\|_{X_{s,b}}$ is a.s. continuous with respect to $T$.
From (\ref{5.07}), we know that $T_{\omega}>0$ a.s. Combining (\ref{5.07}) with the fact that $\|\chi_{[0,T]}\bar{u}\|_{X_{s,b}}$ is $\mathscr{F}_{T}$-measurable, we know that $T_{\omega}$ is a stopping time.
Combining (\ref{5.05}), (\ref{5.06}) with (\ref{5.07}), (\ref{5.08}), we have that $G$ maps the ball of radius 1 in $X_{s,b}^{T_{\omega}}$ into itself and
\begin{eqnarray}
&&\|G(v_{1})-G(v_{2})\|_{X_{s,b}^{T}}\leq
\frac{1}{2}\|v_{1}-v_{2}\|_{X_{s,b}^{T}}
,\label{5.06}
\end{eqnarray}
consequently, $G$ has a unique fixed point, which is the unique process $u$ satisfying (\ref{1.01}) on $[0, T_{\omega}].$
Now we  prove that $u\in C([0,T];H^{s}(\R))$. Since $0<b<\frac{1}{2}$, thus we obtain $\|z(t)\|_{C([0,T];H^{s})}\leq \|z(t)\|_{X_{s,1-b}}$. From Proposition 4.7 of \cite{Richards} and Theorem 6.10 of \cite{PZ}, we know that $\bar{u}\in C([0,T];H^{s}(\R))$. Obviously, we have that
\begin{eqnarray*}
&&\|v\|_{C([0,T];H^{s}}\leq \left\|\frac{1}{3}\int_{0}^{t}U(t-s)\partial_{x}u^{3}ds\right\|_{X_{s,1-b}^{T}}\nonumber\\&&\leq
C\|u\|_{X_{s,b}^{T}}^{3}\leq C(1+\|u_{0}\|_{H^{s}}+C(\omega))^{3}<\infty.
\end{eqnarray*}
Thus, $v\in C([0,T];H^{s}$. In conclusion, we have that $u=z(t)+\bar{u}+v\in C([0,T];H^{s}).$

For  the proof of the  rest of Theorem 1.1,  we refer the readers to Theorem 1.1 of \cite{Bouard-1999,R}.

We have completed the proof of Theorem 1.1.

\bigskip
\noindent {\large\bf 6. Proof of Theorem 1.2 }
\setcounter{equation}{0}

 \setcounter{Theorem}{0}

\setcounter{Lemma}{0}

 \setcounter{section}{6}

In this section, inspired by \cite{R,Richards}, we  prove Theorem 1.2.

Firstly, we consider the following the frequency truncated stochastic PDE

\begin{eqnarray}
\label{6.01}
\begin{cases}
du^{m}(t)=[-\partial_{x}^{3}u^{m}-\frac{1}{3}\partial_{x}((u^{m})^{3})]dt+\Phi_{m} dW(t),\\
u^{m}(x,0)=u_0^{m}(x)=P_{m}u_{0}(x),
\end{cases}
\end{eqnarray}
where $\mathscr{F}_{x}P_{m}u_{0}(x)=\psi\left(\frac{\xi}{m}\right)\mathscr{F}_{x}u_{0}(\xi).$
Obviously, (\ref{6.01})  can be rewritten as follows:
\begin{eqnarray}
u^{m}=U(t)u_{0}^{m}-\frac{1}{3}\int_{0}^{t}S(t-\tau)[(u^{m})^{3}]d\tau+\int_{0}^{t}U(t-\tau)\Phi ^{m}dW(\tau).\label{6.02}
\end{eqnarray}

Firstly, we establish the following Lemmas.

\begin{Lemma}\label{Lemma6.1}Let  $u_{0}(x,\omega)\in  L^{2}
(\Omega; H^{s}(\R)) $ with $s\geq \frac{1}{4}$  and $u_{0}$ be $\mathscr{F}_{0}$  measurable and $\Phi \in L_{2}^{0,\>s}$.
Suppose that $\tilde{\Omega}\subset \Omega$ is such that, for $\omega \in \Omega$, there exists $u^{m}(t)$ which is a solution to (\ref{6.02})
for $t \in [0,T]$ with $T\leq T_{\omega, m}$, where
\begin{eqnarray}
T_{\omega,m}:={\rm \inf} \left\{T>0, 2CT^{1-2b}\left(\|u_{0}^{m}\|_{H^{s}}+2\left\|\psi\left(\frac{t}{T}\right)\bar{u}^{m}\right\|_{X_{s,b}}\right)^{3}\geq 1\right\}.\label{6.03}
\end{eqnarray}
 Then for all $t \in [0,T]$ and any $p\in N$, we have that
 \begin{eqnarray}
 {\rm E}\left(\sup\limits_{t\in [0,T]}\|u^{m}\|_{H_{x}^{1}}^{2p}\chi_{\tilde{\Omega}}\right)\leq C(p,m),\label{6.04}
 \end{eqnarray}
 where $C(p,m)=C\left(p, T, \|u_{0}^{m}\|_{H_{x}^{1}},\|\Phi^{m}\|_{L_{2}^{0,1}}\right)$.
\end{Lemma}
\noindent {\bf Proof.}From Theorem 1.1, we know that there exists a unique solution $u^{m}$ to (\ref{6.01}) for $t\in [0, T_{\omega, N}]$.
Since $T\leq T_{\omega, m}$ inside $\tilde{\Omega}$, we obtain that
\begin{eqnarray}
{\rm E}\left(\sup\limits_{ t\in[0 T]}\|u^{m}\|_{H_{x}^{1}}^{2p}\chi_{\tilde{\Omega}}\right)\leq {\rm E}\left(\sup\limits_{t\in [0, T]}\|u^{m}(t\wedge T_{\omega,m})\|_{H_{x}^{1}}^{2p}\right)\label{6.05}.
\end{eqnarray}
Since $(a+b)^{p}\leq 2^{p-1}(a^{p}+b^{p})$ with $a\geq0,b\geq0,p\geq1$, we have that
\begin{eqnarray}
{\rm E}\left(\sup\limits_{t\in [0, T]}\|u^{m}(t\wedge T_{\omega,m})\|_{H_{x}^{1}}^{2p}\right)\leq\sum_{j=1}^{2}I_{j},\label{6.06}
\end{eqnarray}
where
\begin{eqnarray*}
&&I_{1}=2^{p-1} {\rm E}\left(\sup\limits_{t\in [0, T]}\|u^{m}(t\wedge T_{\omega,m})\|_{L_{x}^{2}}^{2p}\right),\\
&&I_{2}=2^{p-1} {\rm E}\left(\sup\limits_{t\in [0, T]}\|u_{x}^{m}(t\wedge T_{\omega,m})\|_{L_{x}^{2}}^{2p}\right).
\end{eqnarray*}
Obviously,
\begin{eqnarray}
&&I_{2}=2^{p-1} {\rm E}\left(\sup\limits_{t\in [0, T]}(\|u_{x}^{m}(t\wedge T_{\omega,m})\|_{L_{x}^{2}}^{2}-\frac{1}{6}\|u^{m}(t\wedge T_{\omega,m})\|_{L^{4}}^{4}+\frac{1}{6}\|u^{m}(t\wedge T_{\omega,m})\|_{L^{4}}^{4})^{p}\right)\nonumber\\
&&\leq I_{21}+I_{22}.\label{6.07}
\end{eqnarray}
where
\begin{eqnarray*}
&&I_{21}=4^{p-1}{\rm E}\left(\sup\limits_{t\in [0, T]}\left(\|u_{x}^{m}(t\wedge T_{\omega,m})\|_{L^{2}}^{2}-\frac{1}{6}\|u(t\wedge T_{\omega,m})\|_{L^{4}}^{4}\right)^{p}\right),\\
&&I_{22}=\frac{1}{4}\left(\frac{2}{3}\right)^{p}{\rm E}\left(\sup\limits_{t\in [0, T]}\|u(t\wedge T_{\omega,m})\|_{L^{4}}^{4p}\right).
\end{eqnarray*}
By using the interpolation Theorem, we have that
\begin{eqnarray}
I_{22}\leq \frac{2^{p-1}}{4}{\rm E}\left(\sup\limits_{t\in [0, T]}\|u^{m}_{x}(t\wedge T_{\omega,m})\|_{L^{2}}^{2p}\right)+C(p){\rm E}\left(\sup\limits_{t\in [0, T]}\|u(t\wedge T_{\omega,m})\|_{L^{2}}^{6p}\right)\label{6.08}.
\end{eqnarray}
Combining (\ref{6.08}) with (\ref{6.09}), we have that
\begin{eqnarray}
\frac{3}{4}I_{2}\leq I_{21}+C(p){\rm E}\left(\sup\limits_{t\in [0, T]}\|u(t\wedge T_{\omega,m})\|_{L^{2}}^{6p}\right).\label{6.09}
\end{eqnarray}
From (\ref{6.09}), we have that
\begin{eqnarray}
I_{2}\leq \frac{4}{3}I_{21}+C(p){\rm E}\left(\sup\limits_{t\in [0, T]}\|u(t\wedge T_{\omega,m})\|_{L^{2}}^{6p}\right).\label{6.010}
\end{eqnarray}
Combining (\ref{6.06}) with (\ref{6.010}), we have that
\begin{eqnarray}
&&{\rm E}\left(\sup\limits_{t\in [0, T]}\|u^{m}
(t\wedge T_{\omega,m})\|_{H_{x}^{1}}^{2p}\right)\nonumber\\&&
\leq 2^{p-1}{\rm E}\left(\sup\limits_{t\in [0, T]}
\|u(t\wedge T_{\omega,m})\|_{L^{2}}^{2p}\right)\nonumber\\&&\qquad
+\frac{4}{3}I_{21}+C(p){\rm E}\left(\sup\limits_{t\in [0, T]}
\|u(t\wedge T_{\omega,m})\|_{L^{2}}^{6p}\right),\label{6.011}
\end{eqnarray}
Combining (\ref{6.011}) with a  proof similar to (5.3.10) of  Lemma 5.17 of
 \cite{Richards}, we have Lemma 6.1.

We have completed the proof of Lemma 6.1.

\begin{Lemma}\label{Lemma6.2}Let $\alpha=1$ and $u_{0}(x,\omega)\in  L^{2}
(\Omega; H^{s}(\R)) $ with $s\geq \frac{1}{4}$ and $\Phi \in L_{2}^{0,\>s}$
 and $u_{0}$ be $\mathscr{F}_{0}$  measurable.
For any $m$ and any $T_{0}>0$, there exists an almost surely  unique solution
 $u^{m}$  to (\ref{6.02}) for all $t\in [0,T_{0}].$
\end{Lemma}
\noindent {\bf Proof.}
Combining Lemma 6.1 with a proof similar to  Proposition 4.8 of \cite{Richards},
 we have that Lemma 6.2 is valid.

We have completed the proof of Lemma 6.2.

\begin{Lemma}\label{Lemma 6.3}
The sequence $u^{m}$ is bounded in $L^{2}(\Omega, L^{\infty}(0,T_{0};H^{1}(\R))).$ More precisely,
we have that
\begin{eqnarray}
{\rm E}\left(\sup\limits_{t\in [0,T_{0}]}\|u^{m}\|_{H_{x}^{1}}^{2}\right)\leq C\left({\rm E}(\|u_{0}\|_{H^{1}}^{2}), T_{0},\|\Phi\|_{L_{2}^{0,1}}\right).\label{6.012}
\end{eqnarray}
\end{Lemma}
\noindent {\bf Proof.} Let $\mathscr{E}(u^{m})=\|u^{m}\|_{L^{2}}^{6}$.
Applying the It$\hat{o}$ formula to $\mathscr{E}(u^{m})$ yields
\begin{eqnarray}
\|u^{m}\|_{L^{2}}^{6}=\|u_{0}^m\|_{L^{2}}^{6}+6\int_{0}^{t}\|u^{m}\|_{L^{2}}^{4}(u^{m},\Phi^{m}dW)
+\frac{1}{2}\int_{0}^{t}Tr\mathscr{E}^{\prime\prime}(u^{m})(\Phi^{m})(\Phi^{m})^{\star}ds\label{6.013}
\end{eqnarray}
with
\begin{eqnarray*}
\mathscr{E}^{\prime\prime}(u^{m})\phi=24\|u^{m}\|_{L^{2}}^{2}(u^{m},\phi)u^{m}+6\|u^{m}\|_{L^{2}}^{4}\phi.
\end{eqnarray*}
By using a martingale inequality which can be seen in Theorem 3.14 of \cite{PZ}, we have
\begin{eqnarray}
&&{\rm E}\left(\sup \limits_{t\in[0,T_{0}]}\int_{0}^{t}\|u^{m}\|_{L^{2}}^{4}(u^{m},\Phi^{m}dW)\right)\nonumber\\
&&\leq 3{\rm E}\left(\int_{0}^{T_{0}}\|u^{m}\|_{L^{2}}^{8}\|(\Phi^{\star})^{m}u^{m}\|_{L^{2}}ds\right)^{1/2}\nonumber\\
&&\leq \frac{1}{16}{\rm E}\left(\sup\limits_{t\in [0,T_{0}]}\|u^{m}\|_{L^{2}}^{6}\right)+CT_{0}^{3}\|\Phi^{m}\|_{L_{2}^{0,0}}^{6}.\label{6.014}
\end{eqnarray}
By using the definition of trace operator and the Young inequality, we have
\begin{eqnarray}
&&Tr\left(\mathscr{E}^{\prime\prime}(u^{m})\Phi^{m}\Phi^{\star}_{m}\right)\nonumber\\
&&=\sum_{j\in N}\left[24\|u^{m}\|_{L^{2}}^{2}(u^{m},\Phi^{m}e_{j})^{2}+6\|u^{m}\|_{L^{2}}^{4}\|\Phi^{m}e_{j}\|_{L^{2}}^{2}\right]\nonumber\\
&&\leq 30\|u^{m}\|_{L^{2}}^{4}\|\Phi^{m}\|_{L_{2}^{0,0}}^{2}\nonumber\\
&&\leq \frac{1}{12T_{0}}\|u^{m}\|_{L^{2}}^{6}+CT_{0}^{2}\|\Phi^{m}\|_{L_{2}^{0,0}}^{6}.\label{6.015}
\end{eqnarray}
Inserting (\ref{6.014}), (\ref{6.015}) into (\ref{6.013}) yields
\begin{eqnarray}
\|u^{m}\|_{L^{2}}^{6}\leq \|u_{0}^{m}\|_{L^{2}}^{6}+\frac{1}{2}\|u^{m}\|_{L^{2}}^{6}+CT_{0}^{3}\|\Phi^{m}\|_{L_{2}^{0,0}}^{6}\label{6.016}.
\end{eqnarray}
From (\ref{6.016}), we have
\begin{eqnarray}
{\rm E}\left(\sup\limits_{t\in[0,T]}\|u^{m}\|_{L^{2}}^{6}\right)
\leq 2{\rm E}\left(\|u_{0}^{m}\|_{L^{2}}^{6}\right)+CT_{0}^{3}\|\Phi^{m}\|_{L_{2}^{0,0}}^{6}\label{6.017}.
\end{eqnarray}
Let $\mathscr{C}(u^{m})=\|u^{m}\|_{L^{2}}^{8}.$
Applying the It$\hat{o}$ formula to $\mathscr{C}(u_{m})$ yields
\begin{eqnarray}
\|u^{m}\|_{L^{2}}^{8}=\|u_{0}^{m}\|_{L^{2}}^{8}+8\int_{0}^{t}\|u^{m}\|_{L^{2}}^{6}(u^{m},\Phi^{m}dW)
+\frac{1}{2}\int_{0}^{t}Tr\mathscr{C}^{\prime\prime}(u^{m})\Phi^{m}(\Phi^{m})^{\star}ds\label{6.018}
\end{eqnarray}
with
\begin{eqnarray*}
\mathscr{C}^{\prime\prime}(u^{m})\phi=48\|u^{m}\|_{L^{2}}^{4}(u^{m},\phi)u^{m}+8\|u^{m}\|_{L^{2}}^{6}\phi.
\end{eqnarray*}
By using a martingale inequality which can be seen in Theorem 3.14 of \cite{PZ}, we have
\begin{eqnarray}
&&{\rm E}\left(\sup \limits_{t\in[0,T_{0}]}\int_{0}^{t}\|u^{m}\|_{L^{2}}^{6}(u^{m},\Phi^{m}dW)\right)\nonumber\\
&&\leq 3{\rm E}\left(\int_{0}^{T_{0}}\|u^{m}\|_{L^{2}}^{12}\|(\Phi^{\star})^{m}u^{m}\|_{L^{2}}^{2}ds\right)^{1/2}\nonumber\\
&&\leq \frac{1}{16}{\rm E}\left(\sup\limits_{t\in [0,T_{0}]}\|u^{m}\|_{L^{2}}^{6}\right)+CT_{0}^{4}\|\Phi^{m}\|_{L_{2}^{0,0}}^{8}.\label{6.019}
\end{eqnarray}
By using the definition of trace operator and the Young inequality, we have
\begin{eqnarray}
&&Tr\left(\mathscr{C}^{\prime\prime}(u^{m})\Phi^{m}(\Phi^{\star})^{m}\right)\nonumber\\
&&=\sum_{j\in N}\left[48\|u^{m}\|_{L^{2}}^{2}(u^{m},\Phi^{m}e_{j})^{2}+8\|u^{m}\|_{L^{2}}^{4}
\|\Phi^{m}e_{j}\|_{L^{2}}^{2}\right]\nonumber\\
&&\leq 56\|u^{m}\|_{L^{2}}^{6}\|\Phi^{m}\|_{L_{2}^{0,0}}^{2}\nonumber\\
&&\leq \frac{1}{2}\|u^{m}\|_{L^{2}}^{6}+CT_{0}^{3}\|\Phi^{m}\|_{L_{2}^{0,0}}^{8}.\label{6.020}
\end{eqnarray}
Inserting (\ref{6.019}), (\ref{6.020}) into (\ref{6.018}) yields
\begin{eqnarray}
\|u^{m}\|_{L^{2}}^{6}\leq \|u_{0}^{m}\|_{L^{2}}^{6}+\frac{1}{2}\|u^{m}\|_{L^{2}}^{6}+CT_{0}^{3}\|\Phi^{m}\|_{L_{2}^{0,0}}^{6}\label{6.021}.
\end{eqnarray}
From (\ref{6.021}), we have
\begin{eqnarray}
{\rm E}\left(\sup\limits_{t\in[0,T_{0}]}\|u^{m}\|_{L^{2}}^{6}\right)\leq 2{\rm E}\left(\|u_{0}^{m}\|_{L^{2}}^{6}\right)+CT_{0}^{3}\|\Phi^{m}\|_{L_{2}^{0,0}}^{6}\label{6.022}.
\end{eqnarray}
Let
\begin{eqnarray}
H_{2}(u^{m})=\frac{1}{2}\int_{\SR}(u_{x}^{m})^{2}dx-\frac{1}{4}\int_{\SR}(u^{m})^{4}dx\label{6.023}.
\end{eqnarray}
Applying the It$\hat{o}$ formula to $I(u^{m})$  yields
\begin{eqnarray}
&&H_{2}(u^{m})=H_{2}(u_{0}^{m})-\int_{0}^{t}(u_{xx}^{m}+(u^{m})^{3},\Phi dW(s))\nonumber\\&&+\frac{1}{2}\int_{0}^{t}
Tr\left(H_{2}^{\prime\prime}(u^{m})\Phi^{m}(\Phi^{m})^{\star}\right)ds\label{6.024}.
\end{eqnarray}
with
\begin{eqnarray*}
H_{2}^{\prime\prime}(u^{m})\phi=-\phi_{xx}-3(u^{m})^{2}\phi.
\end{eqnarray*}
By using a martingale inequality which can be seen in Theorem 3.14 of \cite{PZ}, we have
\begin{eqnarray*}
&&{\rm E}\left(\sup\limits_{t\in [0,T_{0}]}-\int_{0}^{t}(u_{xx}^{m}+(u^{m})^{3},\Phi^{m} dW(s))\right)\nonumber\\
&&\leq 3{\rm E}\left(\left(\int_{0}^{T_{0}}|(\Phi^{m})^{\star}\left(u_{xx}^{m}+(u^{m})^{3}\right)|^{2}\right)^{1/2}\right).
\end{eqnarray*}
By using the Sobolev embedding $H^{1}\hookrightarrow L^{\infty},$ we have
\begin{eqnarray*}
&&|(\Phi^{m})^{\star}\left(u^{m}_{xx}+(u^{m})^{3}\right)|^{2}=\sum_{j\in N^{+}}\left[(u^{m}_{xx},\Phi^{m}e_{j})+((u^{m})^{3},\Phi^{m}e_{j})\right]^{2}\nonumber\\
&&\leq C\sum_{j\in N^{+}}\left(\|u^{m}\|_{H^{1}}^{2}\|\Phi^{m}e_{j}\|_{H^{1}}^{2}+\|u^{m}\|_{L^{2}}^{4}
\|u^{m}\|_{L^{\infty}}^{2}\|\Phi^{m}e_{j}\|_{L^{\infty}}^{2}\right)\nonumber\\
&&\leq C\left(1+\|u^{m}\|_{L^{2}}^{4}\right)\|u^{m}\|_{H^{1}}^{2}\|\Phi^{m}\|_{L_{2}^{0,1}}^{2}.
\end{eqnarray*}
Consequently, we have
\begin{eqnarray*}
&&{\rm E}\left(\sup\limits_{t\in [0,T_{0}]}-\int_{0}^{t}(u^{m}_{xx}+(u^{m})^{3},\Phi^{m} dW(s))\right)\nonumber\\
&&\leq 3{\rm E}\left(\left(\int_{0}^{T_{0}}|(\Phi^{m})^{\star}\left(u^{m}_{xx}+(u^{m})^{3}\right)|^{2}\right)^{1/2}\right)\nonumber\\
&&\leq C T_{0}^{1/2}\left(1+\|u^{m}\|_{L^{2}}^{2}\right)\|u^{m}\|_{H^{1}}\|\Phi^{m}\|_{L_{2}^{0,1}}\nonumber\\
&&\leq \frac{1}{4}\|u^{m}\|_{H^{1}}^{2}+CT_{0}\|\Phi^{m}\|_{L_{2}^{0,1}}^{2}+
CT_{0}^{2}\|\Phi^{m}\|_{L_{2}^{0,1}}^{4}+C\|u^{m}\|_{L^{2}}^{8}.
\end{eqnarray*}
Thus, by using $H^{1}\hookrightarrow L^{\infty},$ we have
\begin{eqnarray*}
&&Tr\left(H_{2}^{\prime\prime}(u^{m})\Phi^{m}(\Phi^{m})^{\star}\right)\nonumber\\&&=-\sum_{j\in {\rm N}}\int_{\SR}
\left[(\Phi^{m}e_{j})_{xx}\Phi^{m}e_{j}+3(u^{m})^{2}(\Phi^{m}e_{j})^{2}\right]dx\nonumber\\
&&\leq \sum_{j\in {\rm N}}\left(\left|(\Phi^{m}e_{j})_{x}\right|_{L^{2}}^{2}+3\|u^{m}\|_{L^{\infty}}^{2}
\left\|\Phi^{m}e_{j}\right\|_{L^{2}}^{2}\right)
\leq C\|\Phi^{m}\|_{L^{0,1}_{2}}^{2}\left[\|u^{m}\|_{H^{1}}^{2}+1\right]
\end{eqnarray*}
By using  the martingale inequality, we have
\begin{eqnarray*}
{\rm E}\left(\sup\limits_{t\in[0,T_{0}]}-\int_{0}^{t}(u_{xx}^{m}+(u^{m})^{3},\Phi dW(s))\right)\leq 3{\rm E}\left(\left(\int_{0}^{T_{0}}\left|(\Phi^{m})^{\star}\left(u_{xx}^{m}+(u^{m})^{3}\right)\right|^{2}ds\right)^{1/2}\right).
\end{eqnarray*}
Consequently, we have
\begin{eqnarray*}
\frac{1}{2}\int_{0}^{t}Tr\left(H_{2}^{\prime\prime}(u^{m})\Phi^{m}(\Phi^{m})^{\star}\right)ds\leq \frac{1}{4}\|u^{m}\|_{L^{2}}^{4}+CT_{0}^{2}\|\Phi^{m}\|_{L_{2}^{0,1}}^{4}+CT_{0}\|\Phi^{m}\|_{L_{2}^{0,1}}^{2}.
\end{eqnarray*}
Thus, we have
\begin{eqnarray*}
&&{\rm E}\left(\sup\limits_{t\in [0,T_{0}]}H_{2}(u^{m})\right)\nonumber\\&&\leq {\rm E}\left(H_{2}(u_{0}^{m})\right)+\frac{1}{4}{\rm E}\left(\sup \limits_{t\in [0,T_{0}]}\|u_{x}^{m}\|_{L^{2}}^{2}\right)\nonumber\\&&+\frac{1}{4}{\rm E}\left(\sup \limits_{t\in [0,T_{0}]}\|u_{x}^{m}\|_{L^{2}}^{4}\right)+\frac{1}{2}{\rm E}\left(\sup \limits_{t\in [0,T_{0}]}\|u^{m}\|_{L^{2}}^{8}\right)+CT_{0}^{2}\|\Phi^{m}\|_{L_{2}^{0,1}}^{2}+C T_{0}\|\Phi^{m}\|_{L_{2}^{0,1}}^{2}.
\end{eqnarray*}
From the above inequality, by using the interpolation theorem
\begin{eqnarray*}
\|u\|_{L^{2}}^{4}\leq C\|u_{x}\|_{L^{2}}\|u\|_{L^{2}}^{3}+\frac{1}{8}\|u_{x}\|_{L^{2}}^{2}+C\|u\|_{L^{2}}^{6}
\end{eqnarray*}
and (\ref{6.024}), we have
\begin{eqnarray}
&&{\rm E}\left(\sup \limits_{t\in [0,T_{0}]}\|u^{m}_{x}\|_{L^{2}}^{2}\right)\leq
4{\rm E}\left(H_{2}(u_{0}^{m})\right)+C{\rm E}\left(\sup \limits_{t\in [0,T_{0}]}\|u^{m}\|_{L^{2}}^{2}\right)\nonumber\\&&+CT_{0}^{2}\|\Phi^{m}\|_{L_{2}^{0,1}}^{4}
+CT_{0}\|\Phi^{m}\|_{L_{2}^{0,1}}^{2}+C
{\rm E}\left(\sup \limits_{t\in [0,T_{0}]}\|u^{m}\|_{L^{2}}^{4}\right)\nonumber\\&&
\leq C{\rm E}\|u_{0}^{m}\|_{H^{1}}^{2}+C{\rm E}\left(\sup\limits_{t\in [0,T_{0}]}\|u^{m}\|_{L^{2}}^{2}\right)\nonumber\\&&+C{\rm E}\left(\sup\limits_{t\in [0,T_{0}]}\|u^{m}\|_{L^{2}}^{4}\right)+C{\rm E}\left(\sup\limits_{t\in [0,T_{0}]}\|u^{m}\|_{L^{2}}^{6}\right)+C{\rm E}\left(\sup\limits_{t\in [0,T_{0}]}\|u^{m}\|_{L^{2}}^{8}\right)\nonumber\\&&+CT_{0}\|\Phi^{m}\|_{L_{2}^{0,1}}^{2}
+CT_{0}^{2}\|\Phi^{m}\|_{L_{2}^{0,1}}^{4}+CT_{0}^{4}\|\Phi^{m}\|_{L_{2}^{0,1}}^{8}
+CT_{0}^{3}\|\Phi^{m}\|_{L_{2}^{0,1}}^{8}\nonumber\\&&\leq C{\rm E}\left(\|u_{0}^{m}\|_{H^{1}}^{2}\right)+{\rm E}\left(\|u_{0}^{m}\|_{L^{2}}^{6}\right)+
{\rm E}\left(\|u_{0}^{m}\|_{L^{2}}^{8}\right)+C{\rm E}\left(\sup\limits_{t\in[0,T_{0}]}\|u^{m}\|_{L^{2}}^{2}\right)\nonumber\\&&+C{\rm E}\left(\sup\limits_{t\in[0,T_{0}]}\|u^{m}\|_{L^{2}}^{4}\right)+
CT_{0}\|\Phi^{m}\|_{L_{2}^{0,1}}^{2}+CT_{0}^{2}\|\Phi^{m}\|_{L_{2}^{0,1}}^{4}\nonumber\\&&
+C^{6}T_{0}^{3}\|\Phi^{m}\|_{L_{2}^{0,1}}^{6}+
CT_{0}^{4}\|\Phi^{m}\|_{L_{2}^{0,0}}^{8}+CT_{0}^{3}\|\Phi^{m}\|_{L_{2}^{0,0}}^{8}.\label{6.025}
\end{eqnarray}
We define
$\mathscr{D}(u^{m})=\left[\int (u^{m})^{2}dx\right]^{4}.$
Applying the It$\hat{o}$ formula to $\mathscr{D}(u)$ yields
\begin{eqnarray*}
\mathscr{D}(u^{m})=\mathscr{D}(u_{0}^{m})+4\int_{0}^{t}\|u^{m}\|_{L^{2}}^{2}(u^{m},\Phi^{m}dW)+
\frac{1}{2}\int_{0}^{t}Tr\left(\mathscr{D}^{\prime\prime}(u^{m})\Phi^{m}(\Phi^{m})^{\star}\right)ds,
\end{eqnarray*}
where
\begin{eqnarray*}
\mathscr{D}^{\prime\prime}(u^{m})\phi=8(u^{m},\phi)u^{m}+4\|u^{m}\|_{L^{2}}^{2}\phi.
\end{eqnarray*}
By using a computation  similar to (\ref{6.025}), we have
\begin{eqnarray}
{\rm E}\left(\sup\limits_{t\in [0,T_{0}]}\| u^{m}\|_{L^{2}}^{4}\right)\leq 2{\rm E}\left(\| u_{0}^{m}\|_{L^{2}}^{4}\right)+CT_{0}^{2}\|\Phi^{m}\|_{L_{2}^{0,0}}^{4}.\label{6.026}
\end{eqnarray}
In the same way, we have
\begin{eqnarray}
{\rm E}\left(\sup\limits_{t\in [0,T_{0}]}\| u^{m}\|_{L^{2}}^{2}\right)\leq 2{\rm E}\left(\| u_{0}^{m}\|_{L^{2}}^{2}\right)+CT_{0}\|\Phi^{m}\|_{L_{2}^{0,0}}^{2}.\label{6.027}
\end{eqnarray}
Inserting (\ref{6.026}), (\ref{6.027}) into (\ref{6.025}) yields
\begin{eqnarray}
{\rm E}\left(\sup\limits_{t\in [0,T_{0}]}\| u^{m}_{x}\|_{L^{2}}^{2}\right)\leq C{\rm E}\left(\| u_{0}^{m}\|_{H^{1}}+1\right)^{2}+C\left[T_{0}\|\Phi^{m}\|_{L_{2}^{0,1}}^{2}+1\right]^{3}.\label{6.028}
\end{eqnarray}
Combining (\ref{6.027}) with (\ref{6.028}), we have
\begin{eqnarray}
{\rm E}\left(\sup\limits_{t\in [0,T_{0}]}\| u^{m}\|_{H^{1}}^{2}\right)\leq C{\rm E}\left(\| u_{0}^{m}\|_{H^{1}}+1\right)^{2}+C\left[T_{0}\|\Phi^{m}\|_{L_{2}^{0,1}}^{2}+1\right]^{3}.\label{6.029}
\end{eqnarray}
From  (\ref{6.029}), we have that
\begin{eqnarray}
{\rm E} \left(\sup\limits_{t\in [0,T_{0}]}\|u^{m}\|_{H^{1}}^{2}\right)\leq C{\rm E}\left(\| u_{0}\|_{H^{1}}+1\right)^{2}+C\left[T_{0}\|\Phi\|_{L_{2}^{0,1}}^{2}+1\right]^{3}.\label{6.030}
\end{eqnarray}

We have  completed  the proof of  Lemma 6.3.

Now we are in a position to Theorem 1.2.

 From Lemma 6.3,  we know that  after extraction of a subsequence, we can find
a function $\tilde{u} \in L^{2}(\Omega; L^{\infty}(0,T_{0};H^{1}(\R)))$
such that
\begin{eqnarray}
u^{m}\rightharpoonup \tilde{u}\label{6.031}
\end{eqnarray}
in $L^{2}(\Omega; L^{\infty}(0,T_{0};H^{1}(\R)))$ weak star. Moreover,  we have
\begin{eqnarray}
{\rm E} \left(\sup\limits_{t\in [0,T_{0}]}\|\tilde{u}\|_{H^{1}}^{2}\right)\leq C.\label{6.032}
\end{eqnarray}
Let $z^{m}(t)=U(t)u_{0}^{m}$ and $\bar{u}^{m}=\int_{0}^{t}U(t-\tau)\Phi ^{m}d\tau$
and $v^{m}=u^{m}-z^{m}-\bar{u}^{m}$, then for each $m$, $v^{m}$ satisfies the truncated equation
\begin{eqnarray}
v^{m}=\frac{1}{3}\int_{0}^{t}U(t-\tau)\partial_{x}(v^{m}+z^{m}+\bar{u}^{m})^{3}
d\tau=:G_{m}(v^{m}).\label{6.033}
\end{eqnarray}
By repeating the proof of Theorem 1.1,  it is easily checked that $G_{m}$  is a. s. a contraction  on a ball of  radius 1 in $X_{1,\>b}^{\widetilde{T}}$ for any $\widetilde{T}>0$, satisfying
\begin{eqnarray}
2C\widetilde{T}^{1-2b}\left(2+\|u_{0}^{m}\|_{H^{1}}+\|\chi_{t\in [0,\widetilde{T}]}\bar{u}^{m}\|_{X_{1,\> b}}\right) ^{3}\leq 1.\label{6.034}
\end{eqnarray}
Let
\begin{eqnarray*}
D(\omega)=\sup\limits_{0\leq t\leq T_{0}}\|\widetilde{u}\|_{H_{x}^{1}}^{2}.
\end{eqnarray*}
Then
\begin{eqnarray*}
{\rm E}\left(\sup\limits_{0\leq t\leq T_{0}}\|\widetilde{u}\|_{H_{x}^{1}}^{2}\right)\leq C,
\end{eqnarray*}
thus, we derive that $D(\omega)<\infty$ a.s. We consider $\widetilde{T}_{\omega}>0$ satisfying
\begin{eqnarray}
2C\tilde{T_{\omega}}^{1-2b}\left(2+\|u_{0}\|_{H^{1}}+D(\omega)^{1/2}+\|\chi_{t\in [0,\widetilde{T}]}\bar{u}^{m}\|_{X_{1, b}}\right)^{3}\leq 1.\label{6.035}
\end{eqnarray}
Then for any $m$, we have that
\begin{eqnarray*}
\|u_{0}^{m}\|_{H_{x}^{1}}\leq \|u_{0}\|_{H^{1}}
\end{eqnarray*}
and
\begin{eqnarray*}
\left\|\chi_{[0,\widetilde{T}_{\omega}]}\bar{u}^{m}\right\|_{X_{1,b}}\leq \left\|\chi_{t\in [0,\widetilde{T}]}\bar{u}\right\|_{X_{1,b}}.
\end{eqnarray*}
It follows that
(\ref{6.034}) is valid a.s. for any $m$ with 
$\widetilde{T}=\widetilde{T}_{\omega}$. Furthermore, we have that
$\widetilde{T}_{\omega}\leq T_{\omega}$, 
where $T_{\omega}$ is the solution $v$ from Theorem 1.1.
Consequently, $G$ and $G_{m}$ are contractions in $X_{1,b}^{\widetilde{T}_{\omega}}$ for any $m$, where $\widetilde{T}_{\omega}$ satisfies (\ref{6.035}).
Particularly, a unique solution $v\in X_{1,b}^{\widetilde{T}_{\omega}}$
 to (\ref{5.03}) a.s. exists. Moreover, for any $m$, $v^{m}$ and $v$ are
  the unique fixed points of the
contractions $G_{m}$ and $G$, respectively

By using Lemmas 2.4, 2.5, 2.7 and Lemma  4.3, we have that $u^{m}\longrightarrow u$ in $C([0,\widetilde{T}_{\omega}];H^{1}(\mathbf{T}))$ and obtain that
$u=\bar{u}$ for $t\in [0,\widetilde{T}_{\omega}]$ a.s. with the aid of the idea of Section 4.3.2 of \cite{Richards}.  Consequently, we have that
\begin{eqnarray}
\|u(\widetilde{T}_{\omega})\|_{H_{x}^{1}}^{2}\leq \sup \limits_{t\in [0,T_{0}]}\|\bar{u}\|_{H_{x}^{1}}^{2}=D_{\omega}.\label{6.036}
\end{eqnarray}
Combining (\ref{6.035})   with  (\ref{6.036}),  we can construct   a solution on
 $[\widetilde{T}_{\omega},2\widetilde T_{\omega}]$ a.s.
starting from $u(2\widetilde{T}_{\omega})$, we obtain a solution on $[0,T_{0}]$  by reiterating this argument.

We have completed the proof of Theorem 1.2.

\bigskip

\leftline{\large \bf Acknowledgments}

\bigskip

\noindent
We are deeply indebted to the referees for their valuable suggestions which greatly improve
the original version of our paper.


\vskip 10mm
{\large\bf  References}

\begin{thebibliography}{99}

\bibitem{Bouard-1998} A. de Bouard, A. Debussche, On the stochastic Korteweg-de Vries equation, J.
Funct. Anal. 154(1998) 215-251.

\bibitem{Bouard-1999} A. de Bouard, A. Debussche and Y. Tsutsumi, White noise driven Korteweg-de
Vries equation, J. Funct. Anal. 169(1999), 532-558.

\bibitem{Bouard-2007} A. de Bouard, A. Debussche, Random modulation of
 solitons for the stochastic
Korteweg-de Vries equation, {\it Annales de l'Institut Henri Poincar$\acute{e}$
 (C) Analyse Non Lin$\acute{e}$aire,} 24(2007),
251-278.

\bibitem{Bourgain} J. Bourgain, Fourier trandform restriction phenomena
  for certain lattice
subsets and applications to nonlinear evolution equations,
{\it Part II: The KdV equation, Geom. Funct. Anal.}
3(1993), 209-262.

\bibitem{CGG} Y. Chen, H. J.  Gao and B. L. Guo, Well-posedness for
stochastic Camassa-Holm equation,
{\it J. Diff. Eqns.} 253(2012), 2353-2379.

\bibitem{CKST} J. Colliander, M. Keel, G. Staffilani, H. Takaoka, T. Tao,
Global well-posedness for Schr\"odinger equations with derivative,
 {\it SIAM J. Math. Anal.}  33(2001),  649-669.

 \bibitem{CKSTEJDE} J. Colliander, M. Keel, G. Staffilani, H. Takaoka, T. Tao,
 Global well-posedness for KdV in Sobolev spaces of negative index,
   {\it Electr. J. Diff. Eqns.} 26(2001), 1-7.

\bibitem{CKSTMRL}   J.  Colliander, M. Keel, G. Staffilani,
 H. Takaoka, T. Tao,
Almost conservation laws and global rough solutions to a nonlinear
 Schr\"odinger equation,
{\it  Math. Res. Lett.}  9(2002),  659-682.

\bibitem{CKSTJAMS}

J. Colliander, M. Keel, G. Staffilani, H. Takaoka, T. Tao,  Sharp global
well-posedness
 for KdV and modified KdV on $R$  and $T $, {\it  J. Amer. Math. Soc.}
   16(2003),  705-749.


\bibitem{CKS}
J. Colliander, C. E.  Kenig, G. Staffilani, Local well-posedness for
dispersion-generalized Benjamin-Ono equations,
 {\it  Diff. Int. Eqns.}  16(2003),   1441-1472.

\bibitem{GV}
J. Ginibre, G. Velo, Smoothing properties and existence of solutions for
 the generalized Benjamin-Ono equation,  {\it J. Diff. Eqns.}
93(1991), 150-212.


\bibitem{GD}
A. Gr\"unrock, New applications of the Fourier restriction norm method to
well-posedness problems for nonlinear evolution equations,
{\it Dissertation, University of Wuppertal.} 2002.


\bibitem{GIMRN}
A. Gr\"unrock, An improved local well-posedness result for the modified KdV
 equation, {\it  Int. Math. Res. Not.} 61(2004),  3287-3308.

\bibitem{GJMPA}Z. H. Guo,  Global well-posedness of Korteweg-de Vries equation
 in $H^{-3/4} (R),$  {\it J. Math. Pures Appl.}   91(2009),   583-597.

\bibitem{GJDE}Z. H. Guo, Local well-posedness for dispersion generalized
 Benjamin-Ono equations in Sobolev spaces,
 {\it J. Diff. Eqns.}  252(2012),   2053-2084.



\bibitem{Herr} S. Herr, Well-posedness for dispersive equations with derivative
 nonlinearities, {\it Dissertation, Dem Fachbereich  Mathematik  der
University at Dormund vorgelet von.} 2006.


\bibitem{H}
S. Herr,   Well-posedness for equations of Benjamin-Ono type,
 {\it  Illinois J. Math.}  51(2007),   951-976.

\bibitem{HerrCPDE} S. Herr, A. D.  Ionescu, C. E. Kenig, H. Koch,
  A para-differential renormalization technique for nonlinear dispersive equations,
 {\it Comm. Partial Diff. Eqns.}  35(2010),   1827-1875.

\bibitem{Kato} T. Kato, Quasilinear equation of evolution with applications
 to partial differential equations in `` Lect. Notes  in Math.,"
Vol. 448, pp. 27-50, Springer, Berlin, 1975.

\bibitem{KPV1991} C.  Kenig, G.  Ponce, L. Vega, Well-posedness of the initial
 value problem for the Korteweg-de Vries equation,
 {\it J. Amer. Math. Soc.}  4(1991),   323-347.

\bibitem{KPVIUMJ}
 C. Kenig, G.  Ponce, L.  Vega, Oscillatory integrals and regularity of
 dispersive equations, {\it  Indiana Univ. Math. J.}  40(1991), 33-69.

 \bibitem{KPV1993}C.  Kenig, G. Ponce, L.  Vega, Well-posedness and
 scattering results for the generalized Korteweg-de Vries equation via the
  contraction principle, {\it  Comm. Pure Appl. Math.}  46(1993),   527-620.

  \bibitem{KPV1996} C.     Kenig, G. Ponce, L. Vega, A bilinear estimate with
   applications to the KdV equation, {\it  J. Amer. Math. Soc.}  9(1996),  573-603.

\bibitem{KPV2001}C. Kenig, G. Ponce, L. Vega, On the ill-posedness of some
 canonical dispersive equations, {\it  Duke Math. J.}  106(2001),  617-633.

\bibitem{IK}
A. Ionescu, C.  Kenig,  Global well-posedness of the Benjamin-Ono equation
in low-regularity spaces,
{\it  J. Amer. Math. Soc.}  20(2007),  753-798.

\bibitem{KTCPDE} T.  Kappeler, P. Topalov,  Global well-posedness of
 mKdV in $L^{2}(T,R),$ {\it  Comm. Partial Diff. Eqns.}  30(2005),  435-449.

\bibitem{KTD} T.  Kappeler, P.  Topalov,  Global wellposedness of KdV
in $H^{1} (T,R),$  {\it Duke Math. J.}  135(2006),   327-360.

\bibitem{Kis}
 N. Kishimoto, Well-posedness of the Cauchy problem for the Korteweg-de
Vries  equation at the critical regularity, {\it Diff. Int. Eqns.}  22(2009), 447-464.

\bibitem{Koch}H. Koch, N. Tzvetkov,  On the local well-posedness of the Benjamin-Ono
 equation in $H^{s}(R)$, {\it Int. Math. Res. Not.}  26(2003),  1449-1464.

\bibitem{KT}H.  Koch, N. Tzvetkov,  Nonlinear wave interactions for the Benjamin-Ono
 equation, {\it  Int. Math. Res. Not.} 30(2005),   1833-1847.

\bibitem{Molinet} L. Molinet,  Sharp ill-posedness results for the KdV and mKdV
equations on the torus, {\it Advances in Mathematics,}  230(2012), 1895-1930.

\bibitem{MM} L. Molinet, Global well-posedness in the energy space for the Benjamin-Ono
 equation on the circle, {\it  Math. Ann.}  337(2007),  353-383.

 \bibitem{MA} L. Molinet, Global well-posedness in $L^{2}$
  for the periodic Benjamin-Ono equation, {\it  Amer. J. Math.}  130(2008),
    635-683.




\bibitem{MF} L.   Molinet,  Sharp ill-posedness result for the periodic
 Benjamin-Ono equation, {\it  J. Funct. Anal.}  257(2009),   3488-3516.



\bibitem{MP} L.  Molinet, D.  Pilod, The Cauchy problem for the Benjamin-Ono
equation in $L^{2}$    revisited, {\it Anal. PDE}  5(2012),   365-395.




\bibitem{MS} L.  Molinet, J. C.  Saut, N.  Tzvetkov,  Ill-posedness issues
for the Benjamin-Ono and related equations,
 {\it SIAM J. Math. Anal.}  33(2001),   982-988.


 \bibitem{N} T. Nguyen, Power series solution for the modified KdV equations,
 {\it Electron. J. Diff. Eqns.} 15(2008), 359-370.



 \bibitem{NTT}
 K. Nakanishi, H. Takaoka, Y. Tsutsumi, Local well-posedness in low regularity
 of the mKdV equations with periodic bboundary
 condition, {Discrete  Contin. Dyn. Syst.} 28(2010), 1635-1654.

 \bibitem{OR}
 P. Olver, P. Rosenu, Tri-Hamiltonian duality between solitons and solitary-wave
 solutions having compact support,
 {\it Phys. Rev. E.} 53(1996), 1900-1906.


 \bibitem{PZ}
G. Da Prato, J. Zabczyk, Stochastic equations in infinite dimensions, in `` Ency-clopedia
of Mathematics and its applications," Cambridge
University Press, UK, 1992.

\bibitem{R}
 G. Richards, Well-posedness of the stochastic KdV-Burders equation,
  {\it Stochastic Process and their Applications},
 124(2014), 1627-1647.

 \bibitem{Richards} G. Richards, Maximal-in-time behavior of deterministic
  and stochastic dispersive PDEs, Ph.D. Thesis, University of Toronto,
  2012.http://hdl.handle.net/1807/32973.
 \bibitem{SO}
 K. Soonsik, T. Oh, On unconditional well-posedness of modified KdV,
 {\it Int. Math. Res. Not.}
 15(2012), 3509-3534.

 \bibitem{Tao}
 T. Tao, Global well-posedness of the Benjamin-Ono equation in $H^{1}$,
  {\it J. Hyperbolic Differ. Equ.}
 1(2004), 27-49.


 \bibitem{TT}
 H. Takaoka, Y. Tsutsumi, Well-posedness of the Cauchy problem for the modified
 KdV equation with periodic boundary condition,
 {\it Int. Math. Res. Not.}
 56(2004), 3009-3040.


 \bibitem{Ta}
 T. Tao, Multilinear  weighted convolution of $L^{2}$ function and applications
  to nonlinear dispersive equation,
 {\it Amer. J. Math.} 123(2001), 839-908.











\end{thebibliography}
\end{document}